\documentclass[12pt]{article}
\usepackage[utf8]{inputenc}
\usepackage{booktabs}
\usepackage{fullpage}
\usepackage{setspace}
\usepackage[colorlinks,linkcolor=blue, anchorcolor=red,citecolor=blue,CJKbookmarks=true]{hyperref}
\usepackage{epsfig, rotating, amssymb, cite, bbm, amsthm, mathrsfs, dsfont, makeidx,color} %floatrow,
\usepackage{graphicx} % ÓÃÓÚ²åÈëÍ¼ÐÎ
\usepackage{subcaption} % ÓÃÓÚ²åÈë×ÓÍ¼
\usepackage{multirow}
\usepackage{algorithm}
\usepackage{algpseudocode}%
\usepackage{listings}%
\usepackage{bm}%
\usepackage{stmaryrd} %\interleave????????
\usepackage[tbtags]{amsmath}
\allowdisplaybreaks[4]
\usepackage[sectionbib]{natbib}
\newtheorem{theorem}{Theorem}[section]
\newtheorem{corollary}[theorem]{Corollary}

\newcommand{\ve}[1]{{\mbox{\boldmath ${#1}$}}}

\usepackage{moreverb}

\DeclareMathOperator*{\argmin}{arg\,min}

\usepackage{pdflscape}

\begin{document}

%\def\spacingset#1{\renewcommand{\baselinestretch}%
% {#1}\small\normalsize} \spacingset{1}

%\runninghead{X.~Xia.~\textit{ET}~\textit{AL.}}
\title{Biweighted Poisson Subsampling for Convoluted Rank Regression with Massive Data  \footnote{All authors are listed in alphabetical order.}}
% \runtitle{BIPS-DCRR}
\author{ Jialiang Li$^1$
 ~ and ~  Xiaochao Xia$^2$ %\footnote{Corresponding author}  
 ~ and ~ Wei Zhong$^3$ \\  %\thanks{The authors gratefully acknowledge \textit{National Natural Science Foundation of China (No. 11801202)} and \textit{Fundamental Research Funds for the Central Universities (No. 2021CDJQY-047)}.}\\
 \textsl{{ \footnotesize $^{1}$ Department of Statistics and Data Science, National University of Singapore }} \\
 \textsl{{\footnotesize $^{2}$  College of Mathematics and Statistics, Chongqing University} } \\
\textsl{{ \footnotesize $^{3}$ MOE Key Lab of Econometrics, WISE, Department of Statistics and Data Science, Xiamen University }}
}

\date{}
\maketitle

\begin{abstract}
\noindent 

Optimal subsampling efficiently selects the most informative data points, enabling accurate statistical inference while significantly reducing computational burden for massive datasets. However, the existing relevant methods can not directly be applied to pairwise loss problems, particularly for convoluted rank regression (CRR), due to the double summation structure in objective function. To this end, we first propose a new BIweighted Poisson Subsampling (BIPS) framework for such problems through designing a proper weight for a pair of observations instead of for a single observation for objective function. Two concrete inverse probability weighting strategies are considered. Secondly, we focus on the CRR models, under which the BIPS estimator (BIPS-CRR) is formulated. We establish consistency and asymptotic normality for BIPS-CRR, derive its optimal Poisson subsampling probabilities under the L-optimality criterion, and provide a practical algorithm to facilitate implementation. Thirdly, we develop a distributed estimator for CRR that incorporates BIPS as a pilot subsampling strategy. This estimation is globally efficient and is robust to both randomly and non-randomly distributed datasets in distributed computing environments. Extensive simulations and a real-world application demonstrate the excellent finite-sample performance of proposed methodology. Additionally, our BIPS can be readily extended to other U-statistics optimization problems and pairwise learning tasks. 

\end{abstract}

 %\spacingset{1.8} % DON'T change the spacing!

\section{Introduction}
\indent % With the rapid advances in computer and information technology, it is more common for researchers to have access to the data with massive size in various fields such as economics and finance, biology and medicine. However, directly applying traditional statistical methods to these massive datasets faces computational challenges when computational sources are limited. There are at least two solutions to the challenges. One solution is to apply distributed computing methods (\cite{ZDW2013}, \cite{JLY2019}, \cite{Getal2022}). These methods first divide the whole large data set into many relatively small pieces so that traditional methods can be fast applied on each piece of data in parallel, and then aggregate the results from each piece into a final result. Another solution is to utilize subdata selection or subsampling techniques (\cite{YW2021}, \cite{YAY2024}). The subsampling aims to select a small but representative subdata or subsample from the original entire sample data and then existing statistical methods can be implemented on this subsample. %The primary goal of this paper focuses on developing an efficient Poisson subsampling approach and for convoluted rank regression.

Subsampling serves as a powerful tool for handling massive datasets. It performs by utilizing sampling principals to select a relatively small sample from the full sample for fast computation, thus significantly reducing the computational costs. In recent years, plenty of subsampling methods were successfully developed to effectively address the computational challenges imposed by massive data \citep{Cetal2019}. They can be roughly classified to designed-based (model-free) subsampling \citep{Hetal2024, Zetal2023, JV2022, YAY2024} and model-assisted subsampling \citep{MMY2015, WYS2019, WZM2018, WM2021}. The former aims to find representative subdata regardless of the underlying model between the response variable and covariates. The latter is to select most informative samples using specific and trustful models, including leverage value-based subsampling for linear regression \citep{MMY2015}, optimal subsampling with A-optimality for logistic models \citep{WZM2018}, information-based optimal subsampling method \citep{WYS2019} and optimal Poisson subsampling for logistic models \citep{W2018}, among others. While, compared to optimal subsampling with replacement, optimal Poisson subsampling is not only more computationally efficient, due to no need to compute the probabilities for each data point at once, but can also achieve better estimation accuracy \citep{WZW2022}. Important extensions include, but are not limited to, generalized linear models (\citealp{Aetal2021b, SWL2024}), quantile regression (\citealp{Aetal2021a}), asymmetric least squares regression (\citealp{LXZ2024a}), multiplicative regression (\citealp{RZW2023}) and distributed subsampling (\citealp{CMZ2024}). All the optimal subsampling methods as above can solely address single-summation optimizations via inverse probability weighting.

Distributed estimation is another powerful strategy to deal with massive data. Compared to subsampling, it can achieve performance comparable to that of full-sample (global) estimation. In a centralized distributed system, a maasive dataset is distributed on multiple local machines with manageable data sizes, and each can process its own data using traditional statistical models. Only a few quantities are allowed to communicate between central and local machines. A final estimate is obtained by aggregating the local estimates from local machines. There are generally two categories in the literature. One is the one-shot or divide-and-conquer methods \citep{ZDW2013, HH2019, Fetal2019, ZLW2021, LL2023}, which are most highly communication-efficient, but require a large local sample size to ensure global statistical efficiency. The other is iteration methods that need multiple-round communications to trade-off the communication cost and statistical efficiency. Many powerful iteration methods were developed for various models, such as generalized linear model \citep{JLY2019, FGW2023, LMT2025}, quantile regression models \citep{Cetal2020, TBZ2022}, support vector machines \citep{Wetal2019, Xetal2024}, among others. We refer to \cite{Getal2022} for a comprehensive review. The aforementioned distributed approaches all focus on the optimization problems with single-sample loss functions.

%; \citealp{Retal2024}, composite quantile regression (\citealp{SW2022}), ; \citealp{LXZ2024b}

Despite significant progresses in optimal subsampling and distributed estimation, their application to more complex optimization problems, particularly those involving U-statistic-type empirical risk minimization, remains largely unexplored. A notable example is classical rank regression (RR), which relies on a computationally intensive, pairwise loss function \citep{HM2010}. Since RR does not impose distributional assumptions on the response variable, it performs robustly and flexibly, particularly well-suited for data with outliers, heteroscedasticity, or heavy-tailed distribution. Furthermore, compared to quantile regression, RR can achieve greater estimation efficiency for normally distributed data. Several researchers have studied the variable selection problem of under low-dimensional and high-dimensional RR models with relatively small sample sizes \citep{WL2009, L2010, Wetal2020}. Because the RR loss is a non-smooth function, making computation complicated, \citet{ZWZ2024} recently introduced a smoothed variant of RR, termed convoluted rank regression (CRR), through convoluted smoothing techniques. The resulting CRR estimator exhibits improved computational efficiency at the same sample level, and gives better statistical accuracy than RR. % However, RR and CRR relies on pairwise loss functions, which are rather computationally intensive when sample size is large.  

On one hand, directly employing the previous optimal subsampling methods to pairwise loss problems is infeasible due to  three challenges. (i) Pairwise objective functions that take the form of a double summation over all pairs of observations. This structure leads to substantial computational costs, which quickly become prohibitive as the data size increases, especially under memory-constrained environments. (ii) Traditional optimal subsampling only needs to assign weights to individual observations. In contrast, applying inverse probability weighting to the pairwise loss functions requires to design proper biweights for every pair of observations. This increases the difficulty of designing an efficient and principled subsampling strategy. (iii) Due to the U-statistic nature of the pairwise loss function in RR and CRR, deriving optimal subsampling probabilities becomes considerably more difficult than in the models (e.g., quantile regression) with simple additive objectives over independent observations. % For instance, in quantile regression, the objective involves a single summation over i.i.d. data, making the derivation more tractable. 
On the other hand, because subsampling estimators inherently can not achieve the same statistical efficiency to global estimators, how to construct a novel distributed estimation to recover global efficiency under the pairwise loss framework is also an interesting problem, which is largely unexplored. In this paper, to address these challenges and problems, we particularly focus our theoretical development on CRR. Notably, RR can be viewed as a limiting case of CRR as bandwidth tends to zero. 

%\subsection{Main Contribution}
 This paper makes the following major contributions.
\begin{itemize}
  \item Firstly, we propose a new BIweighted Poisson Subsampling (BIPS) framework for pairwise loss problems. Differing from existing methods, our BIPS designs a biweight for a pair of observations instead of for a single observation for objective function. Two concrete weighting strategies are considered via inverse probability weighting. One is with multiplicative weights and the other is with additive weights. The former refers to a subsampling estimator and offers substantially greater computational efficiency while exhibiting slightly lower statistical efficiency than the latter. 

  \item Secondly, under both weighting strategies, we rigorously establish the consistency and asymptotic normality for the BIPS estimator under CRR (BIPS-CRR). We derive optimal subsampling probabilities under the L-optimality criterion that minimizes the trace of an asymptotic covariance matrix of a linear combination of the proposed BIPS-CRR estimator. A practical algorithm for BIPS-CRR is further provided.

  \item Thirdly, under CRR, we propose a distributed estimator that incorporates BIPS as a pilot subsampling to pursuit the global efficiency. Our method is related in spirit to the approach of \cite{JLY2019}, yet it differs in a crucial and novel aspect: while their loss function is constructed solely from data residing on a single local site, our estimator leverages a subsample drawn from the entire dataset via our newly proposed BIPS mechanism. This distinction enables our approach to handle heterogeneous or imbalanced data storage structures, offering flexibility and robustness in distributed computing environments. From a theoretical perspective, we rigorously establish the rate of convergence and prove the asymptotic normality of the estimator, thereby demonstrating its statistical validity and efficiency for large-scale data analysis.

  \item Finally, in comparison to existing methods for optimal Poisson subsampling and distributed estimation, the proposed BIPS addresses broader challenges in both computational feasibility and theoretical development. Notably, our BIPS can be readily extended to other optimization problems involving $U$-statistics, which frequently appear in robust and nonparametric inference (\cite{BC2018}, chapter 2) and in pairwise learning, such as metric learning \citep{JWZ2009, SGH2021}, ranking \citep{CLV2008} and adversarial learning \citep{Wetal2025}. To the best of our knowledge, this is the first study to explore optimal subsampling-based estimation strategies specifically for pairwise $U$-statistic minimization problems. % Our work thus opens a new direction for statistical and machine learning models with pairwise-dependent loss functions.
\end{itemize}

%The rest of this paper is organized as follows. The methodology including the Poisson subsampling  and distributed estimators for CRR are introduced in Section 2. In Section 3, theoretical properties and algorithms are presented. Simulation study is given in Section 4. A real-world application is provided in Section 5. Concluding remarks are given in Section 6. All the proofs of theorems and additional simulation results are relegated to the Appendix.

\section{Methodology}\label{method}

\subsection{Preliminaries}
Suppose we collect the full sample data $\mathcal{D}_n = \{V_i=(Y_i, \mathbf{X}_i)\}_{i=1}^{n}$, consisting of $n$ independent and identically distributed (i.i.d.) observations from the population $V=(Y, \mathbf{X})$ endowed with a probability distribution $\mathbb{P}_{\boldsymbol{\beta}_0}$, where $Y\in \mathbb{R}$ is a response variable and $\mathbf{X}=(X_1,\ldots, X_p)^{\top}\in \mathbb{R}^{p}$ is a $p$-dimensional vector of covariates. 
We assume a linear regression model for the data:
\begin{align}\label{eq:101}
  Y_i = \mathbf{X}_i^{\top} \ve\beta_0 + \epsilon_i, i=1, \ldots, n,
\end{align}
where $\ve\beta_0$ is a $p$-dimensional vector of unknown parameters of interest, and $\epsilon_i$ denotes a random error. We aim to estimate $\ve\beta_0$ in model (\ref{eq:101}) using $\mathcal{D}_n$ when $n$ is extremely large, $p$ is fixed, and $\epsilon_i$ can come from non-normal distribution and may contain outliers. 
In what follows, we are interested in obtaining the estimator of $\ve\beta$ through minimizing a pairwise loss 
\begin{align}\label{eq:101b}
\frac{1}{n(n-1)}\sum_{i=1}^{n}\sum_{i\neq j} \ell(V_i, V_j; \ve\beta)
\end{align}
for robustness considerations, where $\ell(V_i, V_j; \ve\beta)$ represents some loss function that is symmetric in $(V_i, V_j)$. This encompasses two typical examples: classical rank regression (RR) and convoluted rank regression (CRR) estimators that will be considered below. Our approach developed in this section can be applied to other examples associated with the minimization with pairwise loss functions and a related discussion will be given in the last section.

The classical RR estimator of $\ve\beta$ is given by
\begin{align}\label{eq:102}
 \widehat{\ve\beta} =  \argmin_{\boldsymbol{\beta}} \frac{1}{n(n-1)} \sum_{i=1}^{n} \sum_{j\neq i} |Y_{ij}- \mathbf{X}_{ij}^{\top} \ve\beta|,
\end{align}
where $Y_{ij}=Y_i-Y_j$ and $\mathbf{X}_{ij} = \mathbf{X}_i - \mathbf{X}_j, i,j=1,\ldots, n$. It is obvious that the previous loss function becomes $\ell(V_i, V_j; \ve\beta)=|Y_{ij}- \mathbf{X}_{ij}^{\top} \ve\beta|$. This loss function is invariant to a location shift and is well-known to be equivalent to Jaeckel's dispersion function \citep{J1972} with Wilcoxon scores. In this work $\ve\beta_0$ is essentially the minimizer of the population version of the loss function, i.e., $\argmin_{\boldsymbol{\beta}}\mathbb{E}\{ |Y_{ij} -\mathbf{X}_{ij}^{\top} \ve\beta |\}$.
In particular, if $\epsilon_i$s are equipped with distributions symmetric about zero, then $\mathbf{X}_i^{\top} \ve\beta_0$ coincides with the conditional mean. However, our model (\ref{eq:101}) does not need symmetric random error assumption in general. For i.i.d. random errors, $\ve\beta_0 $ can still be interpreted as the effects of the covariates on the conditional mean \citep{Wetal2020}.  % If an intercept term is contained in model (\ref{eq:101}), it can be estimated as $\mathrm{median}(Y_i- \mathbf{X}_i^{\top} \widehat{\ve\beta})$. 
Under some conditions \citep{HM2010}, it has been proved that  $\sqrt{n}(\widehat{\ve\beta}-\ve\beta_0) \stackrel{d}{\to} \mathcal{N}(0, \frac{1}{12\omega^2} \Sigma_{X}^{-1})$, where $\Sigma_{X}=\mathrm{cov}(\mathbf{X})$ denotes the covariance matrix of $\mathbf{X}$, $\omega=\int f(v)^2\mathrm{d}v$, $f(v)$ denotes the density of the random error $\epsilon$ and $\stackrel{d}{\to}$ indicates the convergence in distribution. 
Although the loss is convex, thus ensuring a unique minimizer, the computation of  (\ref{eq:102}) can be challenging due to the non-smoothness and the double summation of the loss. In fact, when the data size is small , one can compute $\widehat{\ve\beta}$ easily using the function {\tt rfit} in R package {\tt Rfit}. However, when the scale of the data is large, this function may not be applicable.

To address the challenge of the non-smoothness of the loss in optimization (\ref{eq:102}), \cite{ZWZ2024} proposed the following CRR estimator: 
\begin{align}\label{eq:501}
 \widehat{\ve\beta}_{h} = \argmin_{\boldsymbol{\beta}} Q_{n, \mathrm{full}}(\ve\beta) := \argmin_{\boldsymbol{\beta}} \frac{1}{n(n-1)} \sum_{i=1}^{n} \sum_{j\neq i} L_h(Y_{ij}- \mathbf{X}_{ij}^{\top} \ve\beta),
\end{align}
where $L_h(\cdot)$ indicates the convoluted rank loss function, defined as the convolution of two functions $L(u)=|u|$ and $K_h(u)$. That is, $L_h(u)=L\circ K_h(u)=\int_{-\infty}^{\infty} |u-v|K_h(v)\mathrm{d}v =u \int_{-u}^{u} K_h(v)\mathrm{d}v - 2 \int_{-\infty}^{u} v K_h(v)\mathrm{d}v $, where $K_h(u)=K(u/h)/h$, $K(u)$ represents a choice of commonly used kernel functions, such as Gaussian kernel $K(u)=\frac{1}{\sqrt{2\pi}}e^{-\frac{u^2}{2}}$ and Epanechnikov kernel $K(u)=0.75(1-u^2)I(|u|<1)$ in nonparametric regression, $h$ represents a bandwidth, and the symbol ``$\circ$" stands for convolution and $I(\cdot)$ denotes the indicator function. Obviously, the previous pairwise loss function reduces to $\ell(V_i, V_j; \ve\beta)= L_h(Y_{ij}- \mathbf{X}_{ij}^{\top} \ve\beta)$. It can be easily verified that the loss $L_h(u)$ is a smooth function with the first derivative $L'_h(u)=2\int_{0}^{u}K_h(v) \mathrm{d} v $ and the second derivative $L''_h(u)=2K_h(u)$.  It is remarkable that if the population parameter is defined as the minimizer of the convoluted rank loss, i.e., $\ve\beta_h^{\ast} = \argmin_{\boldsymbol{\beta}} \mathbb{E} [L_h(Y_{ij} -\mathbf{X}_{ij}^{\top} \ve\beta ) ]$, then it follows that for any $h > 0$, $\ve\beta_h^{\ast}=\ve\beta_0$, as shown in Theorem 1 of \cite{ZWZ2024}. This means that the CRR does not introduce any bias.  Furthermore, \cite{ZWZ2024} established that for any $h>0$, $\sqrt{n}(\widehat{\ve\beta}_{h} - \ve\beta_{0})  \stackrel{d}{\to} \mathcal{N}\left(\mathbf{0}, \frac{\mathbb{E}\{ \mathbb{E}[ L'_h(\epsilon_1 - \epsilon_2)|\epsilon_1]\}^2 }{ [\mathbb{E}L''_h(\epsilon_1 - \epsilon_2)]^2 }\Sigma_{X}^{-1} \right)$. It turns out that in terms of estimation efficiency, the CRR estimator performs arbitrarily close to the RR estimator when $h\to 0$, while behaving asymptotically more efficient than RR for some $h>0$ with a proper kernel choice. Nonetheless, from the computational aspect, both RR and CRR involve a double summation in the loss functions, thus leading to a much higher computing cost. In fact, when applying Newton's iteration method to solve \ref{eq:501}, the entire computation complexity is $O(T(n^2p^2 + p^3))$, where $T$ denotes the total number of iterations. Apparently, when $n$ is huge, the computation becomes infeasible on a personal computer.

\subsection{Biweighted Poisson Subsampling Estimation for CRR/RR}

Subsampling is an effective approach in handling massive data by employing a representative subsample from the entire data set. Compared with uniform subsampling, leverage-based subsampling, and optimal subsampling with replacement, the Poisson subsampling \citep{W2018, WZW2022} has been demonstrated to behave more efficiently. A general procedure of Poisson subsampling proceeds as follows. Start with an empty index set $\mathcal{S}:=\emptyset$. For $i=1,\ldots, n$, the $i$th data point $(\mathbf{X}_i, Y_i)$ is sampled or not according to a Bernoulli experiment with successful probability $\pi_i$: Specifically, we simulate $\delta_i \sim Bernoulli(\pi_i)$. If $\delta_i=1$, $(\mathbf{X}_i, Y_i)$ is sampled and included in the selection index set $\mathcal{S}= \mathcal{S} \cup \{i\}$; otherwise, $\mathcal{S}$ remains the same. Throughout the procedure, the sampling probability $\pi_i$ plays an essential role. We assume $\sum_{i=1}^{n}\pi_i=r$, where $r$ is the expected subsample size. Such a subsampling does not require generating all Bernoulli variables at once, thus saving computational memory for large $n$.

Define a weight $W_i = \delta_i/\pi_i$ for the $i$th data point for $i=1,\ldots, n$. A general Poisson subsampling estimator is the minimizer of an inverse probability weighting loss function, i.e.,
\begin{align}\label{eq:103}
 \widetilde{\ve\beta}  & = \argmin_{\boldsymbol{\beta}} \frac{1}{n} \sum_{i=1}^{n} W_i l(Y_i-\mathbf{X}_i^{\top} \ve\beta)
 = \argmin_{\boldsymbol{\beta}} \frac{1}{n} \sum_{i\in \mathcal{S}} W_i l(Y_i-\mathbf{X}_i^{\top} \ve\beta) ,
\end{align}
where $l(\cdot)$ is a general loss function.  If the squared loss is specified (i.e., $l(u)=u^2$), the resulting estimator, $\widetilde{\ve\beta}$, is referred to the Poisson subsampling estimator for least squares regression \citep{WZW2022}. If the check loss is active (i.e., $l(u)=u[\tau-I(u<0)]$ for some $\tau\in (0,1)$), $\widetilde{\ve\beta}$ is referred to the Poisson subsampling estimator for quantile regression \citep{Aetal2021a}. Furthermore, if the loss is chosen as the two times negative log-likelihood, $\widetilde{\ve\beta}$ is referred to the Poisson subsampling estimator for generalized linear models \citep{WZM2018, WZW2022}.

There are three significant features in (\ref{eq:103}): (i) not all weights are generated at once; (ii) once the weights are specified, the optimization is carried out only on the subsample $\mathcal{S}$, instead of on the full sample; and (iii) the objective function of (\ref{eq:103}) only involves a single summation. Directly incorporating the Poisson subsampling strategy into convoluted rank regression is not straightforward since we have to design a proper weight for each pair of observations with indices $(i,j)$ for $i,j=1,\ldots, n$. 

To address this challenge, we develop the following BIweighted Poisson Subsampling (BIPS) framework. We first describe our BIPS under a general pairwise loss framework and then focus on CRR models. To this end, we assign the biweights $W_{ij}=g(W_i, W_j)$ to the previous pairwise loss function, namely,
\begin{eqnarray}\label{eq:103a}
    Q_{n}(\ve\beta) = \frac{1}{n(n-1)} \sum_{i=1}^{n} \sum_{j\neq i} W_{ij} \ell(V_i. V_j; \ve\beta).
\end{eqnarray}
The bivariate function $g(u, v)$ is prespecified on the support $\{(u,v): u\geq 0, v\geq 0 \}$ such that the biweight $W_{ij}$ satisfies two conditions: (i) $W_{ij}\geq 0$ and (ii) $\mathbb{E}\{ W_{ij}|V_i=(\mathbf{X}_i, Y_i), V_j=(\mathbf{X}_j, Y_j)\}=1$. Condition (i) ensures the non-negativeness of the weighted loss, where letting a part of weights be zero can reduce the computation and nonzero weights correspond to the subsample sampled from the entire data, and condition (ii) is made to guarantee that the conditional expectation of the weighted loss given the original data $\mathcal{D}_n$ equals the original unweighted loss in (\ref{eq:101b}). In fact, there are plenty of choices to meet these two conditions. We consider the multiplicative weights $g(u, v)=uv$ and additive weights $g(u,v)=\frac{u+v}{2}$ in this paper. Once the biweights are determined, one can obtain an estimator of $\ve\beta$ as $\widetilde{\ve\beta}= \argmin_{\boldsymbol{\beta}} Q_{n}(\ve\beta)$. We refer to this procedure as BIPS.

Under CRR, we can formulate the BIPS estimator as
\begin{align}\label{eq:103b}
 \widetilde{\ve\beta}_h =  \argmin_{\boldsymbol{\beta}} Q_{nh}(\ve\beta),
\end{align}
where
%\begin{eqnarray} \label{eq:502b}
$    Q_{nh}(\ve\beta)=\frac{1}{n(n-1)} \sum_{i=1}^{n} \sum_{j\neq i} W_{ij} L_h(Y_{ij}- \mathbf{X}_{ij}^{\top} \ve\beta)$. 
%\end{eqnarray}
% The bivariate function $g(u, v)$ is prespecified on the support $\{(u,v): u\geq 0, v\geq 0 \}$ such that the biweight $W_{ij}$ satisfies two conditions: (i) $W_{ij}\geq 0$ and (ii) $\mathbb{E}\{ W_{ij}|(\mathbf{X}_i, Y_i), (\mathbf{X}_j, Y_j)\}=1$. Condition (i) ensures the non-negativeness of the weighted loss, and condition (ii) is made to guarantee that the conditional expectation of the weighted loss given the original data $\mathcal{D}_n$ equals the original unweighted loss in (\ref{eq:102}). In fact, there are plenty of choices to meet these two conditions. We consider the multiplicative weights $g(u, v)=uv$ and additive weights $g(u,v)=\frac{u+v}{2}$ in this paper. 
With these two choices of weights, we can construct two estimators given by
\begin{align}\label{eq:104}
 \widetilde{\ve\beta}_h^{MW} = \argmin_{\boldsymbol{\beta}} \frac{1}{n(n-1)} \sum_{i=1}^{n} \sum_{j\neq i} W_i W_j L_h(Y_{ij}-\mathbf{X}_{ij}^{\top} \ve\beta)
\end{align}
and
\begin{align}\label{eq:105}
 \widetilde{\ve\beta}_h^{AW} = \argmin_{\boldsymbol{\beta}} \frac{1}{n(n-1)} \sum_{i=1}^{n} \sum_{j\neq i} \frac{W_i+W_j }{2} L_h(Y_{ij}-\mathbf{X}_{ij}^{\top} \ve\beta) 
\end{align}
based on the multiplicative and additive weights, respectively. 

We note that finding $\widetilde{\ve\beta}_h^{MW} $ only needs the subsample $\mathcal{S}=\{i: W_i\neq 0, 1\leq i\leq n\}$, since
\begin{align*}
 \widetilde{\ve\beta}_h^{MW} = \argmin_{\boldsymbol{\beta}} \frac{1}{n(n-1)} \sum_{i\in \mathcal{S}} \sum_{j\in \mathcal{S}, j\neq i} W_i W_j  L_h(Y_{ij}-\mathbf{X}_{ij}^{\top} \ve\beta).
\end{align*}
On the other hand, $\widetilde{\ve\beta}_h^{AW}$ is not a subsampling estimator because $\widetilde{\ve\beta}_h^{AW}$ still requires using the full sample and cannot be obtained from solving $\argmin_{\boldsymbol{\beta}} \frac{1}{n(n-1)} \sum_{i\in \mathcal{S}}\sum_{j\in \mathcal{S}, j\neq i} \frac{1}{2}(W_i+W_j)  L_h(Y_{ij}-\mathbf{X}_{ij}^{\top} \ve\beta)$. We refer to $\widetilde{\ve\beta}_h^{MW} $ as the BIPS estimator for CRR (BIPS-CRR). % We only focus on the practical implementation for BIPS-CRR in the paper.  
From a methodological perspective, both weighted estimators are new and may be of interest. We will still investigate the asymptotic properties of both estimators in the next subsection. In terms of estimation efficiency, one could expect that $ \widetilde{\ve\beta}_h^{AW}$ would be asymptotically more efficient than $\widetilde{\ve\beta}_h^{MW}$ since the former uses all the sample data.

Moreover, when bandwidth approaches zero, the CRR estimators $\widetilde{\ve\beta}_h^{MW}$ and  $\widetilde{\ve\beta}_h^{AW}$ become the corresponding RR estimators, denoted by $ \widetilde{\ve\beta}^{MW}$ and  $\widetilde{\ve\beta}^{AW}$, respectively. We refer to $ \widetilde{\ve\beta}^{MW}$ as the BIPS estimator for RR (denoted as BIPS-RR) which will be implemented in our numerical studies.  

%We remark that we have noted that recently \cite{HX2025} have also developed a similar weighting scheme, however, there are three significant distinctions between their paper and this paper. First, the models considered in both papers are different, where the former is based on RR and ours is mainly for CRR that can regard RR as a special case. Second, subsampling procedures are different, where the former specifically utilizes random perturbation subsampling, while ours focus on Poisson subsampling. In fact, it can be shown that given the same subsampling size, Poisson subsampling can yield a more statistically efficient estimator than that by random perturbation subsampling. Third, technical derivations are quite different. Since our method does involve the derivations for optimal subsampling probabilities, while their method does not, making our entire theories more critically and challenging to establish.  

We remark that our method is non-trivial extension of \cite{HX2025}, who have recently developed a similar weighting scheme. There exists three key distinctions between this and their work. (i) Their work focus on RR model, while this paper is primarily for CRR. (ii) The subsampling procedures differ. Their work specifically employs random perturbation subsampling, whereas ours are with Poisson subsampling. It can be shown that given the same subsampling size, Poisson subsampling yields a more efficient estimator than random perturbation subsampling. (ii) The theoretical derivations are quite distinct. Our work involves deriving optimal subsampling probabilities, while their method does not, making our theoretical results more challenging.

Problem (\ref{eq:103b}) can be solved using the Newton-Raphson method since the loss function is doubly differentiable. 
For such an iterative procedure, the stopping rule can be $\| \widetilde{\ve\beta}_{h}^{(q)}-\widetilde{\ve\beta}_{h}^{(q-1)} \|<\nu $, where $q$ denotes the number of iterations and $\nu$ is a tolerance parameter. We set $\nu=10^{-5}$ in our numerical analysis. Furthermore, the optimal subsampling probabilities $\pi_i$s involved in weights is given in Section 3. Theoretical properties for the estimators and a practical algorithm will be provided.

% Problem (\ref{eq:103b}) can be solved as follows. Let $\widetilde{Y}_{ij}=W_{ij}Y_{ij}$ and $\widetilde{\mathbf{X}}_{ij}=W_{ij}\mathbf{X}_{ij}$, and define a vector $\tilde{\mathbf{y}} =(\widetilde{Y}_{21}, \ldots, \widetilde{Y}_{n1}, Y_{32},\ldots,\widetilde{Y}_{n2}, \ldots, \widetilde{Y}_{n(n-1)})^{\top}\in \mathbb{R}^{n(n-1)/2}$ and a matrix $\tilde{\mathbf{x}} \in \mathbb{R}^{n(n-1)/2 \times p}$, similarly. Thus, problem (\ref{eq:103b}) can be reformulated as $\argmin_{\boldsymbol{\beta}} \frac{2}{n(n-1)} \sum_{i=1}^{n(n-1)/2} |\tilde{y}_{i}- \tilde{\mathbf{x}}_{i}^{\top}\ve\beta|$, where $\tilde{y}_{i}$ and $\mathbf{\mathbf{x}}_{i}^{\top}$ correspond the $i$th element of $\tilde{\mathbf{y}}$ and the $i$th row of $\tilde{\mathbf{x}}$, respectively. This becomes a least absolute deviation (LAD) problem, which can be easily solved using the R package {\tt quantreg}.
% As previously illustrated, addressing the LAD problem with large dataset could be exceedingly time-consuming and potentially prohibitive. Fortunately, this situation is unlikely to occur because most elements of $\tilde{\mathbf{y}}$ and many rows of $\tilde{\mathbf{y}}$ are exactly zero due to many weights being zero, as long as the expected subsample size is small relative to the full sample size.

\subsection{Distributed Estimation with Biweighted Poisson Subsampling }
%Although the proposed Poisson subsampling method is computationally efficient when $n$ is exceedingly large, the estimation efficiency cannot achieve the global efficiency.
From our theoretical development, we will show that the convergence rate of our BIPS-CRR estimator is $r^{-1/2}$, slower than that of the global CRR estimator. Thus we further propose a distributed estimation for CRR with the above biweighted Poisson subsampling to pursuit the global estimation efficiency.

Suppose the entire dataset $\mathcal{D}_n$ is divided into $M$ subsets of equal size and stored on $M$ machines. The $m$th machine possesses the subdata $\mathcal{D}_{(m)} = \{ (\mathbf{X}_{i}, Y_{i}), i\in \mathcal{I}_m\}$, where $\mathcal{I}_m$ denotes the set of indices and its cardinality is $|\mathcal{I}_m|=n_m = n/M$ for $m=1,\ldots, M$. We propose a distributed estimation procedure in which a master machine communicates with $M$ local machines. Information exchange is permitted between the master and each local machine; however, communication among local machines is not allowed.

Firstly, according to the BIPS-CRR estimation given in the previous section, we can obtain the subsampling estimator $\widetilde{\ve\beta}_{h}$ as an initial estimator, i.e., $\widehat{\ve\beta}^{(0)}: = \widetilde{\ve\beta}_{h}$. Next, based on the distributed data, we compute $\{\nabla Q_{n,m}(\widehat{\ve\beta}^{(0)})\}_{m=1}^{M}$, which are local gradients of
\begin{align}\label{eq:509}
    Q_{n,m}(\ve\beta) =  \frac{1}{n_m(n_m-1)} \sum_{i\in \mathcal{I}_m} \sum_{ j\in \mathcal{I}_m, j\neq i} L_h(Y_{ij}-\mathbf{X}_{ij}^{\top}\ve\beta),
\end{align}
evaluated at $\widehat{\ve\beta}^{(0)}$. Such calculation can be performed in parallel on local machines. We then aggregate these gradients by $\nabla \widetilde{Q}_n(\widehat{\ve\beta}^{(0)}): = \frac{1}{M}\sum_{m=1}^{M} \nabla Q_{n,m}(\widehat{\ve\beta}^{(0)})$ on the master machine. We now define a surrogate loss
\begin{align}\label{eq:510}
\widetilde{Q}_n(\ve\beta) := Q_{nh}(\ve\beta) - [\nabla Q_{nh}(\widehat{\ve\beta}^{(0)}) - \nabla \widetilde{Q}_n(\widehat{\ve\beta}^{(0)})]^{\top} \ve\beta,
\end{align}
where $Q_{nh}(\ve\beta)$ is given in (\ref{eq:103b}) based on the Poisson subsample with weights $\{W_i, i\in \mathcal{S}\}$, and $\nabla Q_{nh}(\ve\beta)$ denotes its gradient. Using a Taylor series expansion we can easily see that the surrogate loss could closely approximate the loss function based on the full sample. Hence, our updated estimator is simply the minimizer of the surrogate loss $\widetilde{Q}_n(\ve\beta) $, i.e.,
\begin{align}\label{eq:511}
\widehat{\ve\beta}^{(1)} =  \argmin_{\boldsymbol{\beta}} ~\widetilde{Q}_{n}(\ve\beta).
\end{align}
This can be solved by a Newton-Raphson procedure  on the master machine. Such a nonlinear optimization is performed with only the Poisson subsamples. We refer to the estimator $\widehat{\ve\beta}^{(1)}$ as the distributed CRR estimator with biweighted Poisson subsampling (BIPS-DCRR) in this paper. It is noted that compared to the CSL estimator in \cite{JLY2019}, our estimator has two key distinctions: (i) $\nabla \widetilde{Q}_n(\widehat{\ve\beta}^{(0)}) \neq \nabla  Q_{n, \mathrm{full}}(\widehat{\ve\beta}^{(0)})$ and (ii) the solved optimization uses the subsamples by Poisson subsampling rather than local samples.

\section{Theoretical Properties}\label{theory}

In this section, we present the theoretical results of the proposed estimators in Section \ref{method}. 
In our theoretical analysis, the expected size, $r$, of the subsample is allowed to increase with the full sample size, $n$, i.e., $r\to\infty$ as $n\to\infty$. Note that the condition of $r/n \to 0$ is not required. Instead, we allow $r/n \to c_0\in [0,1]$ as $n\to\infty$.

\subsection{Theoretical Properties for BIPS-CRR}
In this subsection, we establish the asymptotic properties of the proposed estimators in (\ref{eq:104}) and (\ref{eq:105}). We need the following regularity conditions.

% Denote the probability density function (PDF) and the cumulative distribution function (CDF) of the random error $\epsilon$ in model (\ref{eq:101}) by $f(\cdot)$ and $F(\cdot)$, respectively. 
%% (A5) -> (A4), (A6)-> (A5), (A7)->(A6)

\begin{itemize}
  \item[(A1)] Let $\Theta$ be a compact parameter space, a convex subset of $\mathbb{R}^{p}$. The true parameter $\ve\beta_{0}$ lies in the interior of $\Theta$, i.e., $\ve\beta_{0} \in \mathrm{int}(\Theta)$, and $R:=\sup_{\boldsymbol{\beta} \in \Theta} \|\ve\beta - \ve\beta_0 \| >0$

  \item[(A2)] There exist universal positive constants $C_l, l=1,\ldots, 6,$ such that
     (i) $\mathbb{E}\{ \|\mathbf{X}\|^3 \}\leq C_1$,
     (ii) $\mathbb{E}\{ \frac{r}{n^2} \sum_{i=1}^{n} \frac{1}{\pi_i} \} \leq C_2$,
     (iii) $\mathbb{E}\{ \frac{r}{n^2} \sum_{i=1}^{n} \frac{1}{\pi_i} \|\mathbf{X}_i\|^2 \} \leq C_3$,
     (iv) $\mathbb{E}\{ \frac{r^2}{n^3} \sum_{i=1}^{n} \frac{1}{\pi_i^2} \|\mathbf{X}_i\|^2 \} \leq C_4$,
     (v) $\mathbb{E}\{ \frac{r}{n^2} \sum_{i=1}^{n} \frac{1}{\pi_i} \|\mathbf{X}_i\|^3 \} \leq C_5$,
     and (vi) $\mathbb{E}\{ \frac{r^2}{n^3} \sum_{i=1}^{n} \frac{1}{\pi_i^2} \} \leq C_6$,

  \item[(A3)] (i) There exist two positive constants $M_1$ and $M_2$ such that $\sup_{v} f(v) \leq M_1$ and $\sup_{v} |f'(v)|\leq M_2$. (ii) Furthermore, $\epsilon$ and $\mathbf{X}$ are independent.

   \item[(A4)]  The kernel function $K(u)$ satisfies (i) $K(-u)=K(u) ~\forall u\in \mathbb{R}$; (ii) $\kappa_{u} = \sup_{u\in \mathbb{R}} K(u)<\infty$ and $\exists \varsigma_0>0, s.t.~ \kappa_{l} \equiv \inf_{u\in [-\varsigma_0, \varsigma_0]} K(u) >0$; (iii) $\int_{-\infty}^{\infty} K(u)\mathrm{d}u =1$; and (iv) $\exists \alpha_0\in (0,1], \tilde{L}_0>0, ~s.t. ~|K(u)-K(v)|\leq \tilde{L}_0|u-v|^{\alpha_0}, \forall u,v\in \mathbb{R}$.

  \item[(A5)] Let $\bar{f}(\cdot)$ be the probability density function of $\epsilon_{12}=\epsilon_1-\epsilon_2$. Assume (i) $|\bar{f}(u)-\bar{f}(v)| \leq \tilde{L}_1 |u-v|^{\alpha_1}, ~\forall x,y\in \mathbb{R}$, for some constants $\alpha_1\in (0,1]$ and $\tilde{L}_1>0$; and (ii) $\mu_{u} = \sup_{u\in \mathbb{R}} \bar{f}(u)<\infty$ and $\exists \varsigma_1>0, s.t.~ \mu_{l} \equiv \inf_{u\in [-\varsigma_1, \varsigma_1]} \bar{f}(u) >0$.

  \item[(A6)] Assume that $\Omega_{\pi h} := \mathrm{plim}_{n\to \infty} \frac{r}{n^2} \sum_{i=1}^{n} \frac{1}{\pi_i} \left[\mathbb{E}\{L'_h(\epsilon-\epsilon')|\epsilon=\epsilon_i\}\right]^2 (\mathbf{X}_i -\mathbb{E}\mathbf{X})(\mathbf{X}_i -\mathbb{E}\mathbf{X})^{\top} $ is a positive definite matrix, where $\epsilon'$ is an independent copy of $\epsilon$, and $\mathrm{plim}$ indicates convergence in probability.
\end{itemize}

 Assumptions (A1)-(A6) are mild. Assumption (A1) is a regularity condition on the parametric space $\Theta$ and the true parameter $\ve\beta_0$ (\cite{Yetal2022} and \cite{JLY2019}). Assumption (A2) is a technical requirement on the moments of the covariates and the sampling probabilities for data points. %It is similar to the assumption in classic sampling techniques (see, e.g., Breidt and Opsomer 2000; Berger and Torres 2016).
It can protect the optimization from being dominated by data points with extremely small subsampling probabilities. %Furthermore, Assumption (A2) is satisfied if $P\{c_1 r/n \leq \pi_i \leq c_2 r/n, \forall i=1,\ldots, n\}\to 1$ holds for two constants $0<c_1<c_2$ as $(n,r)\to\infty$. 
In Assumption (A3), the first part requires the density of $\epsilon$ and its first derivative to be bounded from above. This is similar to Assumption 1 in \cite{ZWZ2024}. The second part restricts the model error and covariates to be independent. This is a standard condition in rank regression literature \citep{ZWZ2024, L2010}.

%% Remark 5 -> Remark 2,  Remark 6 -> Remark 3,  Remark 7 -> Remark 5, Remark 2 -> Remark 4
%% Theorem 1 -> Corollary 1, Theorem 2 -> Corollary 2,
  Assumption (A4) is a standard requirement on kernel function $K(\cdot)$. Common choices in kernel regression, such as the Gaussian kernel, $K(u)= (2\pi)^{-1/2} \exp(-u^2/2)$, and the Epanechnikov kernel, $K(u)=0.75 (1-u^2)I(|u|\leq 1)$, can satisfy these conditions.  Assumption (A5) concerns the density of the difference of two i.i.d. model errors. Assumptions (A4) and (A5) are the same as those in \cite{ZWZ2024}. Assumption (A6) is to ensure the positive definiteness of the asymptotic covariance matrix of the proposed estimators. In particular, if uniform subsampling probabilities (i.e., $\pi_i = \frac{r}{n}, i=1,\ldots, n$) are used, $\Omega_{\pi,h}$ in Assumption (A6) becomes $\Omega_{\pi h} = \mathbb{E}\{ \mathbb{E}[ L'_h(\epsilon_1-\epsilon_2)|\epsilon_1]^2 \} \Sigma_{X}$.

\begin{theorem} \label{th04}
If Assumptions (A1)-(A6) are satisfied, we have (i) (rate of convergence) $ \widetilde{\ve\beta}_{h}^{MW} - \ve\beta_{0} = O_p(r^{-1/2}) $; and
(ii) (asymptotic normality)
\begin{eqnarray*}
\sqrt{r} (\widetilde{\ve\beta}_{h}^{MW} - \ve\beta_{0}) \stackrel{d}{\to} \mathcal{N}\left(\mathbf{0}, \frac{\Sigma_{X}^{-1} \Omega_{\pi h} \Sigma_{X}^{-1}}{ [\mathbb{E}L''_h(\epsilon_1 - \epsilon_2)]^2 } \right).
\end{eqnarray*}
\end{theorem}

\begin{theorem} \label{th05}
If Assumptions (A1)-(A6) are satisfied, we have  (i) (rate of convergence) $ \widetilde{\ve\beta}_{h}^{AW} - \ve\beta_{0} = O_p(r^{-1/2}) $; and
(ii) (asymptotic normality)
\begin{eqnarray*}
\sqrt{r} (\widetilde{\ve\beta}_{h}^{AW} - \ve\beta_{0}) \stackrel{d}{\to} \mathcal{N}\left(\mathbf{0}, \frac{\Sigma_{X}^{-1}(3c_0 \mathbb{E}\{ \mathbb{E}[ L'_h(\epsilon_1-\epsilon_2)|\epsilon_1]^2 \} \Sigma_{X} + \Omega_{\pi h})\Sigma_{X}^{-1}}{4[\mathbb{E}L''_h(\epsilon_1 - \epsilon_2)]^2}\right).
\end{eqnarray*}
\end{theorem}

 It can be verified that as $h\to 0$, the limits of $\mathbb{E}L''_h(\epsilon_1 - \epsilon_2)$ and $\mathbb{E}\{ \mathbb{E}[ L'_h(\epsilon_1-\epsilon_2)|\epsilon_1]^2 \}$ in asymptotic variances of Theorems \ref{th04} and \ref{th05} approach $2\omega$ and $1/3$, respectively, and that $\mathbb{E}[ L'_h(\epsilon_1-\epsilon_2)|\epsilon_1 = u] \to 2 F(u)-1$ as $h\to 0$ for any fixed $u\in \mathbb{R}$, where $F(\cdot)$ stands for the cumulative distribution function of the random error $\epsilon$. In fact, we have the following two corollaries under rank regression.

\begin{corollary} \label{th01}
Suppose that $\Omega_{\pi} := \mathrm{plim}_{n\to\infty}  \frac{r}{n^2}\sum_{i=1}^{n} \frac{1}{\pi_i} (\mathbf{X}_i-\mathbb{E} \mathbf{X}_i)(\mathbf{X}_i-\mathbb{E} \mathbf{X}_i)^{\top}[F(\epsilon_i) -\frac{1}{2}]^2 $ is positive definite. Under Assumptions (A1)-(A3), we have (i) (rate of convergence) $ \widetilde{\ve\beta}^{MW} - \ve\beta_{0} = O_p(r^{-1/2}) $; and (ii) (asymptotic normality)
\begin{eqnarray*}
\sqrt{r} (\widetilde{\ve\beta}^{MW} - \ve\beta_{0}) \stackrel{d}{\to} \mathcal{N}\left(\mathbf{0}, \frac{1}{\omega^{2}}\Sigma_{X}^{-1} \Omega_{\pi} \Sigma_{X}^{-1} \right).
\end{eqnarray*}
\end{corollary}

\begin{corollary} \label{th02}
Under the assumptions of Corollary \ref{th01}, we have (i) (rate of convergence) $ \widetilde{\ve\beta}^{AW} - \ve\beta_{0} = O_p(r^{-1/2}) $; and
(ii) (asymptotic normality)
\begin{eqnarray*}
\sqrt{r} (\widetilde{\ve\beta}^{AW} - \ve\beta_{0}) \stackrel{d}{\to} \mathcal{N}\left(\mathbf{0}, \frac{1}{4\omega^{2}}\Sigma_{X}^{-1}(\frac{c_0}{4}\Sigma_{X} + \Omega_{\pi})\Sigma_{X}^{-1}\right).
\end{eqnarray*}
\end{corollary}

% This implies that Theorems \ref{th01} and \ref{th02} can be viewed a special case of Theorems \ref{th04} and \ref{th05}, respectively, when the bandwidth is very small and close to zero.

 From Corollaries \ref{th01} and \ref{th02}, we can make two observations. First, when $c_0=0$ (i.e., $r/n\to 0$ as $n\to\infty$) holds, the asymptotic variance of $\widetilde{\ve\beta}^{AW}$ equals $\frac{1}{4\omega^{2}}\Sigma_{X}^{-1} \Omega_{\pi} \Sigma_{X}^{-1}$, which is only a quarter of the asymptotic variance of $\widetilde{\ve\beta}^{MW}$ no matter how we specify the Poisson subsampling probabilities. Using the notion of the asymptotic relative efficiency (ARE), we have $\mathrm{ARE}(\widehat{\ve\beta}^{AW}, \widehat{\ve\beta}^{MW})=4$. Second, if uniform subsampling probabilities are used, then $\sqrt{r} (\widetilde{\ve\beta}^{MW} - \ve\beta_{0}) \stackrel{d}{\to} \mathcal{N}\left(\mathbf{0}, \frac{1}{12 \omega^2} \Sigma_{X}^{-1}\right)$ and $\sqrt{r} (\widetilde{\ve\beta}^{AW} - \ve\beta_{0}) \stackrel{d}{\to} \mathcal{N}\left(\mathbf{0}, \frac{3c_0+1}{48 \omega^2} \Sigma_{X}^{-1}\right)$.
This leads to $\mathrm{ARE}(\widehat{\ve\beta}^{AW}, \widehat{\ve\beta}^{MW})= \frac{4}{3 c_0+1}$ with $c_0 \in [0,1]$. Further, in the extreme case where $c_0=1$ (i.e., $r=n$ and $\pi_i=1$), we have $\Omega_{\pi}= \frac{1}{12}\Sigma_{X}$, implying that $\widetilde{\ve\beta}^{MW}$ and $\widetilde{\ve\beta}^{AW}$ have the same asymptotic covariance matrix, $\frac{1}{12\omega^{2}}\Sigma_{X}^{-1}$, agreeing with the classical results for rank regression \citep{HM2010}.

  In the special case of uniform subsampling probabilities, we can obtain
\begin{align} \label{eq:503}
   \mathrm{ARE}(\widetilde{\ve\beta}_{h}^{AW}, \widetilde{\ve\beta}^{AW}) =  \mathrm{ARE}(\widetilde{\ve\beta}_{h}^{MW}, \widetilde{\ve\beta}^{MW}) =  \frac{[\mathbb{E}L''_h(\epsilon_1 - \epsilon_2)]^2}{ 12 \omega^2 \mathbb{E}\{ \mathbb{E}[ L'_h(\epsilon_1-\epsilon_2)|\epsilon_1]^2 \}},
\end{align}
and
\begin{align} \label{eq:504}
  \mathrm{ARE}(\widetilde{\ve\beta}_{h}^{AW}, \widetilde{\ve\beta}_h^{MW}) = \frac{4}{3c_0+1}.
\end{align}
On one hand, we have $\mathrm{ARE}(\widetilde{\ve\beta}_{h}^{AW}, \widetilde{\ve\beta}_h^{MW})>1$ as long as $c_0\neq 1$, implying that $\widetilde{\ve\beta}_{h}^{AW}$ is more efficient than $\widetilde{\ve\beta}_{h}^{MW}$ for uniform subsampling with $r<n$. On the other hand, it can be shown that $\mathrm{ARE}(\widetilde{\ve\beta}_{h}^{MW}, \widetilde{\ve\beta}^{MW}) \to 1$ and $\mathrm{ARE}(\widetilde{\ve\beta}_{h}^{AW}, \widetilde{\ve\beta}^{AW}) \to 1$ as $h\to 0$. Since zero may not be the optimal choice for $h$, from Theorem 3 of \cite{ZWZ2024}, the right-hand side of (\ref{eq:503}) equals $1+ C_{K,f} h^2 +o(h^2)$ in the neighborhood of $0$, where $C_{K,f}$ is a constant that merely relies on the kernel $K(\cdot)$ and the density $f(\cdot)$ of the error $\epsilon$. The constant $C_{k,f}$ can be positive using Epanechnikov kernel and normal distribution of the random error. This indicates that both
$\mathrm{ARE}(\widetilde{\ve\beta}_{h}^{AW}, \widetilde{\ve\beta}^{AW})$ and $\mathrm{ARE}(\widetilde{\ve\beta}_{h}^{MW}, \widetilde{\ve\beta}^{MW})$ can be strictly greater than one for sufficiently small $h$. This further indicates that BIPS-CRR can be more efficient than the subsampling estimator for classic RR.

Next, we provide the optimal subsampling probabilities $\pi_i$s involved in BIPS-CRR. From Theorems \ref{th04} and \ref{th05}, one can see that the asymptotic covariance matrices, $\mathrm{cov}(\widetilde{\ve\beta}_h^{MW})$ and $\mathrm{cov}(\widetilde{\ve\beta}_h^{AW})$ both depend on subsampling probabilities $\ve\pi :=\{\pi_1,\ldots, \pi_n\}$ via $\Omega_{\pi}$. To determine an optimal $\ve\pi$, one can minimize the asymptotic mean squared error (AMSE) of the two estimators with respect to $\ve\pi$. That is equivalent to minimizing the traces of $\mathrm{cov}(\widetilde{\ve\beta}_h^{MW})$ and $\mathrm{cov}(\widetilde{\ve\beta}_h^{AW})$. Such a method is usually referred to as A-optimality criterion. Another method is to adopt the L-optimality criterion, that is, the optimal $\ve\pi$ is achieved by minimizing the trace of the asymptotic covariance matrix of $\Sigma_{X}\widetilde{\ve\beta}_h^{MW}$, i.e., minimizing $\mathrm{tr}(\Omega_{\pi h})$. 
%some linear transformation of the parameter estimator, say $L^{\top} \widetilde{\ve\beta}_h^{MW}$, where $L$ is a proper matrix. This amounts to minimizing $\mathrm{tr}(L\Sigma_X^{-1}\Omega_{\pi h} \Sigma_{X}^{-1} L^{\top})$. If $L=\Sigma_{X}$, the L-optimality criterion deduces to
Compared to the L-optimality, the A-optimality criterion may require more complicated computation to obtain the estimator of $\Sigma_X$. The following result gives the optimal subsampling probabilities for the proposed estimators in terms of L-optimality criterion.

\begin{theorem} \label{th06}
Define $\hbar_i^{\dag} := \|\mathbf{X}_i-\mathbb{E} \mathbf{X} \||\mathbb{E}\{L'_h(\epsilon_i-\epsilon')|\epsilon_i\}|$, $i=1,\ldots,n$.
For the estimators $\widetilde{\ve\beta}_{h}^{MW}$ and $\widetilde{\ve\beta}_{h}^{AW}$, the L-optimal subsampling probabilities that minimize $\mathrm{tr}(\Omega_{\pi h,n})$, where
\begin{align}  \label{eq:505}
 \Omega_{\pi h, n} =  \frac{r}{n^2} \sum_{i=1}^{n} \frac{1}{\pi_i} \left[\mathbb{E}\{L'_h(\epsilon-\epsilon')|\epsilon=\epsilon_i\}\right]^2 (\mathbf{X}_i -\mathbb{E}\mathbf{X})(\mathbf{X}_i -\mathbb{E}\mathbf{X})^{\top},
\end{align}
are
\begin{eqnarray} \label{eq:506}
 \pi_{h,n,i}^{L, opt} = r \frac{ \hbar_i^{\dag} \wedge H  }{\sum_{j=1}^{n} \hbar_j^{\dag}  \wedge H }, i=1,\ldots, n,
\end{eqnarray}
where $a \wedge b = \min(a,b)$, and $ H= \sum_{i=1}^{g} \hbar^{\dag}_{(i)}/(g+r-n)$, in which $\hbar^{\dag}_{(1)}\leq \cdots\leq \hbar^{\dag}_{(n)}$ are the order statistics of $\{\hbar_i^{\dag}\}_{i=1}^{n}$, and $g$ is an integer such that
\begin{align}
  \frac{\hbar^{\dag}_{(g)}}{\sum_{i=1}^{g} \hbar^{\dag}_{(i)}} < \frac{1}{g+r-n},  %  \label{eq:506a} \\
  ~\text{and}~ \frac{\hbar^{\dag}_{(g+1)}}{\sum_{i=1}^{g+1} \hbar^{\dag}_{(i)}} \geq \frac{1}{ g+r-n+1},  \label{eq:506b}
\end{align}
where we define $\hbar^{\dag}_{(n+1)} = \infty$.
\end{theorem}
Theorem \ref{th06} provides the optimal Poisson subsampling probabilities under the L-optimality criterion. In fact, this result can be generalized to a setting where one can minimize the trace of the asymptotic covariance of the linear combination of the proposed estimator, i.e., $\min_{\boldsymbol{\pi}} \mathrm{tr}(\mathrm{cov}(L\widetilde{\ve\beta}_h^{MW}))$ for a proper matrix $L$. It contains the A-optimality and L-optimality as two special cases. In addition, we provide a practical algorithm to implement the optimal Poisson subsampling probabilities for BIPS-CRR and the details can be found in Section \ref{algorithm}.

\subsection{Theoretical Properties for Distributed Estimation}

% condition (A8) -> (A7); Remark 8 -> Remark 6; Remark 9 -> Remark 7; 
To establish the asymptotic properties of the estimator $\widehat{\ve\beta}^{(1)}$ in (\ref{eq:511}), we need an additional condition (A7) given below.  

\begin{itemize}
  % \item[(A7)] Assume that $\mathbb{E} [L''_h(\epsilon_{12})]^8 <\infty$, $\sup_{1\leq j\leq p} \mathbb{E} X_{j}^{16} <\infty$ and $\mathbb{E}[\frac{r^8}{n^8} \sum_{i=1}^{n} \frac{1}{\pi_i^8} ] <\infty $.
  % \item[(A7)] Assume that $\mathbb{E} [L''_h(\epsilon_{12})]^4 <\infty$, $\sup_{1\leq j\leq p} \mathbb{E} X_{j}^{16} <\infty$ and $\mathbb{E}[\frac{r^4}{n^4} \sum_{i=1}^{n} \frac{1}{\pi_i^4} ] <\infty $. 
\item[(A7)] (i) $\sup_{1\leq j\leq p} \mathbb{E} X_{j}^{4} <\infty$ and (ii) $\max_{1\leq n} \frac{1}{n\pi_i} = O_p(r^{-1})$.  
\end{itemize}
Condition (A7)(i) requires the fourth-order moments of covariates bounded, which is much weaker than that of \cite{CXZ2026} (their Remark 1), the Assumption PD in \cite{JLY2019} and condition (C3) in \cite{LXZ2024a}. Condition (A7)(ii) is an assumption on subsampling probabilities, which is commonly used in subsampling literature, see Assumption 5 in \cite{Yetal2022} for example.

\begin{theorem} \label{th07}
Suppose that Assumptions (A3)(ii) and (A4)-(A7) hold, and the initial estimator satisfies $\|\widehat{\ve\beta}^{(0)} - \ve\beta^{\ast}\|=O_p(b_n)$, where $b_n$ is a positive sequence that tends to zero as $n\to\infty$. We have (i) (rate of convergence)
\begin{align*}
  \| \widehat{\ve\beta}^{(1)} - \ve\beta_{0} \| = O_p\left( \frac{1}{\sqrt{n}} + (\frac{1}{r}+\frac{1}{\sqrt{n}})b_n + b_n^{1+\alpha_0}  \right),  % + \frac{\sqrt{M}}{n}
\end{align*}
where $\alpha_0$ is a constant given in condition (A4); and (ii) (asymptotic normality) if $\sqrt{n}b_n \max( \frac{1}{r}, b_n^{\alpha_0})=o(1)$ holds as $(n,r)\to\infty$, then
\begin{align*}
\sqrt{n}(\widehat{\ve\beta}^{(1)} - \ve\beta_{0})  \stackrel{d}{\to} \mathcal{N}\left(\mathbf{0}, \frac{\mathbb{E}\{ \mathbb{E}[ L'_h(\epsilon_1 - \epsilon_2)|\epsilon_1]\}^2 }{ [\mathbb{E}L''_h(\epsilon_1 - \epsilon_2)]^2 }\Sigma_{X}^{-1} \right).
\end{align*}
\end{theorem}

Part (i) suggests that the proposed distributed estimator, $\widehat{\ve\beta}^{(1)}$, can improve the convergence rate of the initial estimator, $\widehat{\ve\beta}^{(0)}$, with a contraction factor of $\frac{1}{r} + \frac{1}{\sqrt{n}} + b_n^{\alpha_0}$. When the Poisson subsampling estimator $\widehat{\ve\beta}_h$ is used as the initial estimator, as described by Algorithm \ref{alg03} in Section \ref{algorithm}, it follows that $b_n=r^{-1/2}$ and the convergence rate in part (i) of Theorem \ref{th07} is $O_p(\frac{1}{\sqrt{n}}+ \frac{1}{r^{(1+\alpha_0)/2}})$. To achieve a global convergence rate, this amounts to requiring that $r\geq n^{1/(1+\alpha_0)}$. Furthermore, if our subsampling estimator is chosen as the initial estimator, the number of machines $M$ does not affect the convergence rate in the asymptotic sense. However, if the first local estimator is used as the initial estimator (i.e. $b_n=\sqrt{\frac{M}{n}}$), $M$ could affect the convergence rate.  Part (ii) indicates that the proposed estimator, $\widehat{\ve\beta}^{(1)}$, can behave as well as the global estimator under the assumption $\sqrt{n}b_n \max( \frac{1}{r}, b_n^{\alpha_0})=o(1)$. When $b_n=r^{-1/2}$, this assumption reduces to that $r$ is strictly greater than $n^{1/(1+\alpha_0)}$ in order. The asymptotic distribution result of part (ii) is exactly the same as that of the global estimator proposed in \cite{ZWZ2024}. It is worth noting that the proof of this theorem employs a different technique from \cite{JLY2019}. A key technical challenge in our proof is to handle the term $\nabla \widetilde{Q}_{n}(\widehat{\ve\beta}^{(0)}) - \nabla Q_{n, \text{full}} (\widehat{\ve\beta}^{(0)}) $, which is zero in the framework of \cite{JLY2019}, but non-zero in our pairwise loss problem. All the proofs of the above results are given in Appendix C.2 of supplementary file. 

% \smallskip
% \noindent{\textit{Remark 7}} Algorithm \ref{alg03} is an one-step algorithm as it solely performs with one round of communication between the master and local machines, and hence, is highly communication-efficient. The global rate of convergence can be guaranteed by requiring $r\geq n^{1/(1+\alpha_0)}$. However, if this is not met, the property may fail. There are two manners to remain this property. One manner is to increase the size of subsample and the other is to perform multiple rounds of communications. For the latter, one needs to perform more than $\lfloor \frac{\log n}{ \alpha_0 \log r} \rfloor$ rounds of communications between the master and local machines, where $\lfloor a\rfloor$ represents the largest integer not exceeding the number $a$.

\section{Practical Implementation}\label{algorithm}

\subsection{The Algorithm for BIPS-CRR}

Note that there are unknown quantities $F(u)$, $\mathbb{E}\{L'_h(\epsilon_{i}-\epsilon)|\epsilon_i\}$ and $\mu_{X}=\mathbb{E} \mathbf{X}$ involved in the optimal probabilities in Theorem \ref{th06}.  We consider a pilot subsampling approach to approximate the optimal subsampling probabilities. We aim to select a pilot subsample of a relatively small size $r_0$ from the original data using a deterministic subsampling method.  The uniform replacement subsampling is usually recommended for ease of implementation.  Let $\mathcal{P}$ be the set of indices of the pilot subsample and  $\mathcal{D}_{r_0}=\{ (\mathbf{X}_i,Y_i),i\in \mathcal{P} \}$ be the resulting pilot subsample. We can obtain the pilot estimator of $\mu_{X} \triangleq \mathbb{E}\mathbf{X}$ as $\widehat{\mu}_{X,r_0} = \frac{1}{r_0} \sum_{i\in \mathcal{P}} \mathbf{X}_i$ and the pilot subsampling CRR estimator of $\ve\beta$ as
\begin{align} \label{eq:401}
   \widehat{\ve\beta}_{h,r_0} = \argmin_{\boldsymbol{\beta}} \frac{1}{r_0(r_0-1)} \sum_{i\in \mathcal{P}}\sum_{j\in \mathcal{P}, j\neq i} L_h(Y_{ij}-\mathbf{X}_{ij}^{\top}\ve\beta ).
\end{align}
Next, we output the residuals $\{\hat{\epsilon}_{i, r_0} = Y_i - \mathbf{X}_i^{\top} \widehat{\ve\beta}_{r_0}, i\in\mathcal{P}\}$, and obtain the pilot estimator of $F(u)$ as $\widehat{F}_{r_0}(u)=\frac{1}{r_0}\sum_{i\in \mathcal{P}}I(\hat{\epsilon}_{i, r_0}\leq u)$. It remains to obtain a pilot estimator of $J_h(u) \triangleq  \mathbb{E}\{L'_h(\epsilon-\epsilon')|\epsilon=u\}$. Note that $J_h(u) = 2\int_{-\infty}^{\infty} F(u- ht)K(t)\mathrm{d}t -1$. A natural estimator is to plug $\widehat{F}_{r_0}(\cdot)$ into the integral $J_h(u)$ to obtain the pilot estimator of $J_h(u)$. However, numerically computing this integral is very inaccurate and unstable. Instead, we utilize the average of $L'_h(u-\hat{\epsilon}_{i,r_0})$ based on the pilot residuals $\{\hat{\epsilon}_{i,r_0}=Y_i-\mathbf{X}_i^{\top} \widehat{\boldsymbol{\beta}}_{i,r_0}, i=1,\ldots, r_0\}$ as the estimator $\widehat{J}_h(u)$, i.e., $\widehat{J}_h(u)=\frac{1}{r_0} \sum_{i=1}^{r_0} L'_h(u-\hat{\epsilon}_{i,r_0})$. 
As a result, a feasible estimator of $ \pi_{h,n,i}^{L, opt} $ can be constructed as
\begin{eqnarray} \label{eq:507}
 \widehat{\pi}_{h, n,i}^{L,opt} = r \frac{ \|\mathbf{X}_i-\widehat{\mu}_{X,r_0} \| | \widehat{J}_{h, r_0}(\widehat{\epsilon}_i) | \wedge \widehat{H}  }{\sum_{j=1}^{n} \|\mathbf{X}_j-\widehat{\mu}_{X,r_0} \| |\widehat{J}_{h, r_0}(\widehat{\epsilon}_j)| \wedge \widehat{H} }, i=1,\ldots, n,
\end{eqnarray}
where $\widehat{H}$ corresponds to the empirical version of $H$ given in Theorem \ref{th06} where the three unknown quantities are substituted by their pilot estimators as above. In terms of the A-optimality, we can similarly obtain a feasible optimal Poisson subsampling probability as
\begin{eqnarray} \label{eq:508}
 \widehat{\pi}_{h, n,i}^{A,opt} = r \frac{ \|\widehat{\Sigma}_{X,r_0}^{-1}(\mathbf{X}_i-\widehat{\mu}_{X,r_0}) \| | \widehat{J}_{h, r_0}(\widehat{\epsilon}_i) | \wedge \widehat{H}  }{\sum_{j=1}^{n} \|\widehat{\Sigma}_{X,r_0}^{-1}(\mathbf{X}_j-\widehat{\mu}_{X,r_0}) \| |\widehat{J}_{h, r_0}(\widehat{\epsilon}_j)| \wedge \widehat{H} }, i=1,\ldots, n,
\end{eqnarray}
where $\widehat{\Sigma}_{X,r_0}$ is the sample covariance of $\mathbf{X}$ based on $\mathcal{D}_{r_0}$ and $\widehat{H}$ can be obtained similarly. To accelerate computation, the denominators in (\ref{eq:507}) and (\ref{eq:508}) can be replaced by $n \widehat{\Psi}_1$ and $n\widehat{\Psi}_2$, respectively, where $ \widehat{\Psi}_1$ ($\widehat{\Psi}_2$) is the average of $\{\|\mathbf{X}_j-\widehat{\mu}_{X,r_0} \| |\widehat{J}_{h, r_0}(\widehat{\epsilon}_j)| \wedge \widehat{H}, i\in \mathcal{P}\}$ ($\{\|\widehat{\Sigma}_{X,r_0}^{-1}(\mathbf{X}_j-\widehat{\mu}_{X,r_0}) \| |\widehat{J}_{h, r_0}(\widehat{\epsilon}_j)| \wedge \widehat{H}, i\in \mathcal{P}\}$). We summarize the procedure in Algorithm \ref{alg02}.

\begin{algorithm}
\caption{A practical algorithm for BIPS-CRR.}\label{alg02}
\begin{algorithmic}[1]
\State \textbf{Input} $r_0, r, \varrho$ and $\mathcal{D}_n$.
\State \textbf{Step 1: Pilot subsampling}: Initialization: ${\mathcal{D}}_{r_0}=\emptyset; $
\For {$i=0,\dots,n$ }
  \State Generate a Bernoulli random variable $\delta_i \sim Bernoulli(\pi_i)$ with $\pi_i =r_0/n$;			
  \If {$\delta_i=1$ }
	\State Update ${\mathcal{D}}_{r_0}={\mathcal{D}}_{r_0} \cup\{(y_i,\mathbf{X}_i)\}; $			
  \EndIf
\EndFor
\State Compute estimators  $\widehat{\ve\beta}_{h,r_0}$, $\widehat{F}_{r_0}(\cdot)$, $\widehat{\mu}_{X,r_0}$ and $\widehat{\Sigma}_{X,r_0}^{-1}$ based on the pilot sample ${\mathcal{D}}_{r_0}$.  \\ 	
\State \textbf{Step 2: Approximated optimal Poisson subsampling}:  Initialization: ${\mathcal{D}}_{r}=\emptyset; $
\For{$i=1,\dots,n$}
  \State  Calculate the optimal subsampling probability $\pi_{h, n,i}^{\ast, opt} = \widetilde{\pi}_{h, n,i}^{opt} \wedge 1$, where
  \begin{align} \label{eq:405b}
   \widetilde{\pi}_{h, n,i}^{opt} = (1-\varrho)\widehat{\pi}_{h, n,i}^{opt} + \varrho \frac{r}{n},
  \end{align}
  and $\widehat{\pi}_{h, n,i}^{opt}$ can be either (\ref{eq:507}) or (\ref{eq:508}).
  \State Generate a Bernoulli random variable $\delta_i \sim Bernoulli(\pi_{h, n,i}^{\ast, opt})$;
  \If {$\delta_i=1$ }
     \State Update ${\mathcal{D}}_{r}={\mathcal{D}}_{r}\cup\{(y_i,\bm{x}_i)\}; $
     \State Record $W_i = 1/\pi_{h, n,i}^{\ast, opt}$;	
   \EndIf
\EndFor
\State \textbf{Step 3: Estimation}: Let $\mathcal{S}$ be the set of indices of ${\mathcal{D}}_{r}$  and find the minimizer, $\widetilde{\ve\beta}_{h,r}$, of
\begin{align*}
   \min_{\boldsymbol{\beta}} \frac{1}{n(n-1)} \sum_{i\in \mathcal{S}}\sum_{j\in \mathcal{S}, j\neq i} W_iW_j L_h(Y_{ij}-\mathbf{X}_{ij}^{\top}\ve\beta).
\end{align*}
\State \textbf{Output}: $\widetilde{\ve{\beta}}_{h} = \widetilde{\ve\beta}_{h,r}$.
\end{algorithmic}
\end{algorithm}

In the supplementary file, we provide more discussion and intuition for this algorithm in Section A1 and give a practical algorithm for BIPS-RR in Appendix A.3. 

\subsection{The Algorithm for BIPS-DCRR}

We summarize the previous distributed estimation procedure stated in Section 2.3 in Algorithm \ref{alg03}. We remark that the initial estimator used in Algorithm \ref{alg03} can be replaced by another different initial estimator, such as the local estimator computed from the data on the first machine in \cite{JLY2019}.

\begin{algorithm}
% \SetAlgoNoLine  %???????
\caption{A practical algorithm for BIPS-DCRR.}\label{alg03}
\begin{algorithmic}[1]
\State \textbf{Input} $r_0$, $r$, $\varrho$ and $\mathcal{D}_n$.
\State \textbf{Step 1: Poisson Subsampling}: Draw a Poisson subsample $\mathcal{D}_{r}$ via Algorithm \ref{alg02} and store it on the master machine.
\State Compute an initial estimator $\ve\beta^{(0)} := \widetilde{\boldsymbol{\beta}}_{h}$ based on $\mathcal{D}_{r}$ and transmit it to local machines.
\State \textbf{Step 2: Distributed Estimation}: Divide the entire dataset $\mathcal{D}_n$ evenly into $M$ subsets $\{\mathcal{D}_{(m)}\}_{m=1}^{M}$ stored on $M$ machines.
\For {$m=1,\ldots, M$}
  \State Compute the gradient $\nabla Q_{n,m}(\ve\beta^{(0)})$, where $Q_{n, m} (\boldsymbol{\beta})$ is defined in (\ref{eq:509}).
  \State Transmit it to the master machine.
\EndFor
\State On the master machine, compute the aggregated gradient $\nabla \widetilde{Q}_n(\ve\beta^{(0)}) = M^{-1}\sum_{m=1}^{M} \nabla Q_{n,m}(\ve\beta^{(0)}) $ and $ \nabla Q_{nh}(\boldsymbol{\beta}^{(0)})$.
\State Find $\widehat{\ve\beta}^{(1)} = \argmin_{\boldsymbol{\beta}} ~\widetilde{Q}_{n}(\ve\beta)$, where $\widetilde{Q}_{n}(\ve\beta)$ is defined in (\ref{eq:510}).
%\begin{align*}
%   \min_{\boldsymbol{\beta}}\left\{ Q_{nh}(\boldsymbol{\beta}) - \left[\nabla Q_{nh}(\ve\beta^{(0)}) - \nabla \widetilde{Q}_{n}(\ve\beta^{(0)})  \right]^{\top} \boldsymbol{\beta} \right\}.
%\end{align*}\\
\State \textbf{Output} $\widehat{\ve\beta}^{(1)}$.
\end{algorithmic}
\end{algorithm}

Algorithm \ref{alg03} is an one-step algorithm as it solely performs with one round of communication between the master and local machines, and hence, is highly communication-efficient. To ensure the global rate of convergence, $r\geq n^{1/(1+\alpha_0)}$ is required in theory. There are two practical approaches to satisfy this. One way is to increase the size of subsample and the other is to perform multiple rounds of communications. For the latter, one needs to perform more than $\lfloor \frac{\log n}{ \alpha_0 \log r} \rfloor$ rounds of communications between the master and local machines, where $\lfloor a\rfloor$ represents the largest integer not exceeding the number $a$. An additional algorithm based on multiple rounds of communications is provided in Appendix B.3 of supplementary file.

\section{Simulation Study}
We conduct a simulation study to illustrate the finite-sample performance of our proposed estimation approaches in previous sections. We mainly evaluate the performance of estimators in terms of their mean squared errors (MSE). Throughout this section, we fix $n=50,000$, $r_0=200$ and number of repetitions $N=200$, and use Epanechnikov kernel in CRR.

\textbf{Example 1} We generate a full dataset of $n$ independent observations $\mathcal{D}_n=\{(y_i, \mathbf{X}_i)\}_{i=1}^{n}$ from the model: $y_i = \mathbf{X}_i^{\top} \boldsymbol{\beta}_0 + \epsilon_i,i=1,\ldots, n$, where $\mathbf{X}_i$ is simulated from a multivariate normal distribution $\mathcal{N}(\mathbf{0}, \Sigma)$. Here, we consider two structures for the covariance matrix $\Sigma=(\sigma_{lj})_{1\leq l,j\leq p}$: (i)AR(0.5): autoregression correlation $\sigma_{lj} = 0.5^{|l-j|}$; and (ii) CS(0.5): compound symmetric correlation $\sigma_{lj}=0.5$ if $l\neq j$ and $\sigma_{lj}=1$ if $l=j$.
We consider three cases for the error distribution: (i) the standard normal distribution $\epsilon_i \sim \mathcal{N}(0,1)$; (ii) the student's $t$ distribution with four degrees of freedom $\epsilon_i \sim \sqrt{2}t(4)$; and (iii) the mixed normal distribution $\epsilon_i \sim 0.95\mathcal{N}(0,1)+0.05\mathcal{N}(0, 50)$. We set the true parameter vector as $\boldsymbol{\beta}_0=(\sqrt{3}/2,\sqrt{3}/2, \sqrt{3}/2, \sqrt{3}/2,\sqrt{3}/2, \sqrt{3}/2)^{\top}$.

Due to the page limit, we put all the simulation results and analysis for our proposed BIPS-CRR in section B1 of the supplementary file, along a detailed comparison with other four subsample methods: (BIPS-RR) the Poisson subsampling estimator under RR, (PS-QR) the Poisson subsampling method under quantile regression \citep{Aetal2021a}, (PS-LM) the Poisson subsampling under linear models \citep{WZW2022}, and (IBOSS) the information-based optimal subdata selection \citep{WYS2019}. The simulation results demonstrate that our BIPS-CRR method consistently outperforms BIPS-RR, with both approaches showing superior performance compared to existing alternatives.

In this section we further investigate the finite-sample performance of the proposed distributed estimator BIPS-DCRR based on Algorithm \ref{alg03}. We compare our BIPS-DCRR with three popular distributed estimations: (CSL) the communication-efficient distributed estimation based on a surrogate loss \citep{JLY2019}, (DLSA) the distributed estimation based on least squares approximation \citep{ZLW2021}, and (AE) the averaged estimation that aggregates the estimators from local machines via simple average. Note that the CSL and the proposed BIPS-DCRR methods can be implemented with multiple rounds of communications, but the DLSA and AE methods need only one round of communication, which are also known as one-shot methods. For the CSL method, we use the first local estimator as the initial estimator of $\ve\beta$ as in \cite{JLY2019}, while for our BIPS-DCRR method, we use the proposed BIPS-CRR estimator with the L-optimality criterion as an initial estimator. We consider the following two distributed storage schemes for the entire dataset, and under each scheme, we compare the performance of the aforementioned four distributed estimation methods in terms of MSE.

\begin{itemize}
  \item (Scheme 1) We randomly partition the full dataset $\mathcal{D}_n$ generated by Example 1 into $M$ blocks with approximately equal size, and then each block is stored on a local machine.
  \item (Scheme 2) Following Example 1, we first generate $\mathcal{D}_n$ and compute $\{U_i=\sum_{j=1}^{p} X_{ij},i=1,\ldots, n\}$ as well as their ordered statistics $\{U_{(1)}\leq \cdots \leq U_{(n)}\}$. We then allocate the data points $\{(\mathbf{X}_i, Y_i), i\in \mathcal{I}_m\}$ with $\mathcal{I}_m = \{i: U_{((m-1)n/M+1)} \leq U_{(i)} \leq U_{(mn/M)}\}$ to the $m$th local machine for $m=1,\ldots, M$.
\end{itemize}

\begin{figure} %[htp]
\centering
    \includegraphics[height=8cm, width=12cm]{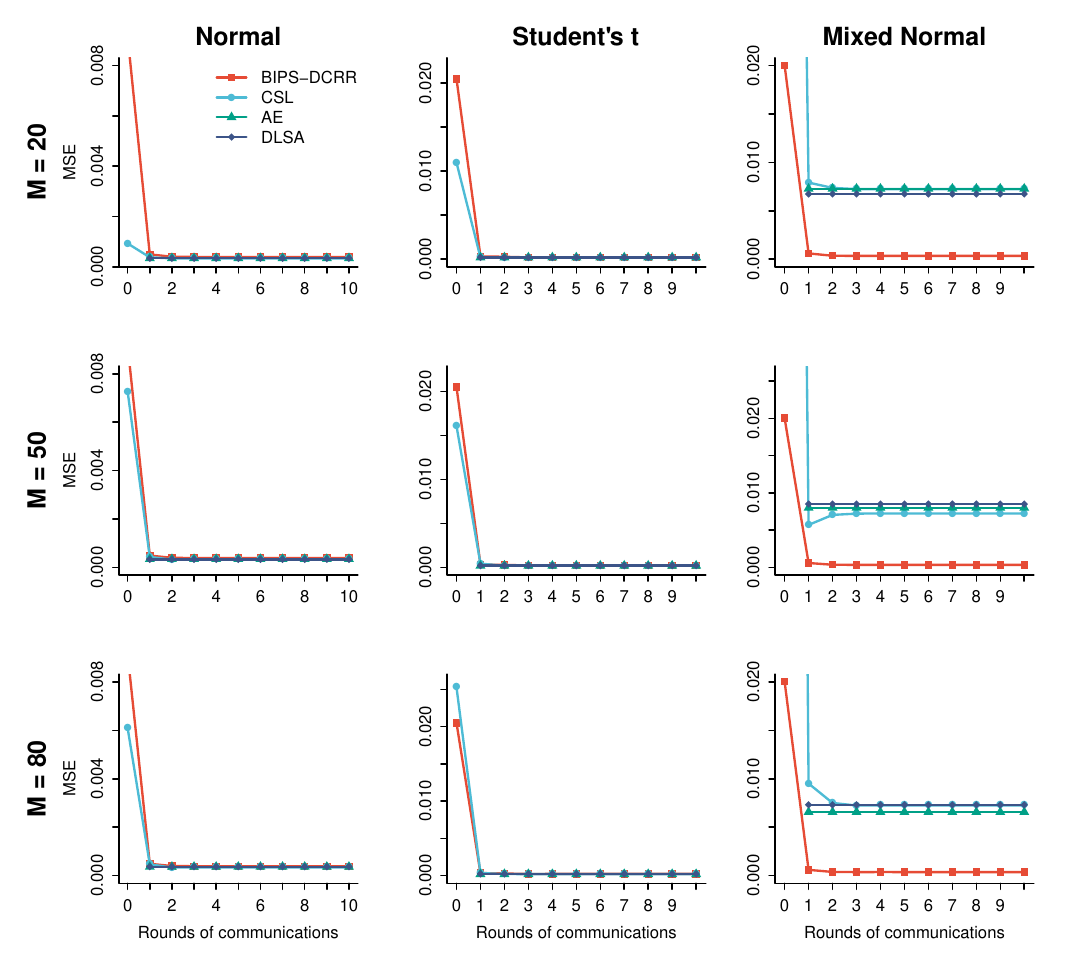} \\
\caption{Results of MSEs for four distributed estimation methods under Scheme 1, where we fix $n=50,000, r_0=200, r=1,000$, $\varrho=0.35$, $h=1$ and $\mathbf{X}\sim \mathrm{AR}(0.5)$.  (Normal): $\epsilon\sim \mathcal{N}(0,1)$, (Student's t): $\epsilon\sim \sqrt{2}t(4)$, and (Mixed Normal): $\epsilon\sim 0.95\mathcal{N}(0,1)+0.05\mathcal{N}(0,50)$. }
\label{fig:505}
\end{figure}

\begin{figure} %[htp]
\centering
    \includegraphics[height=8cm, width=12cm]{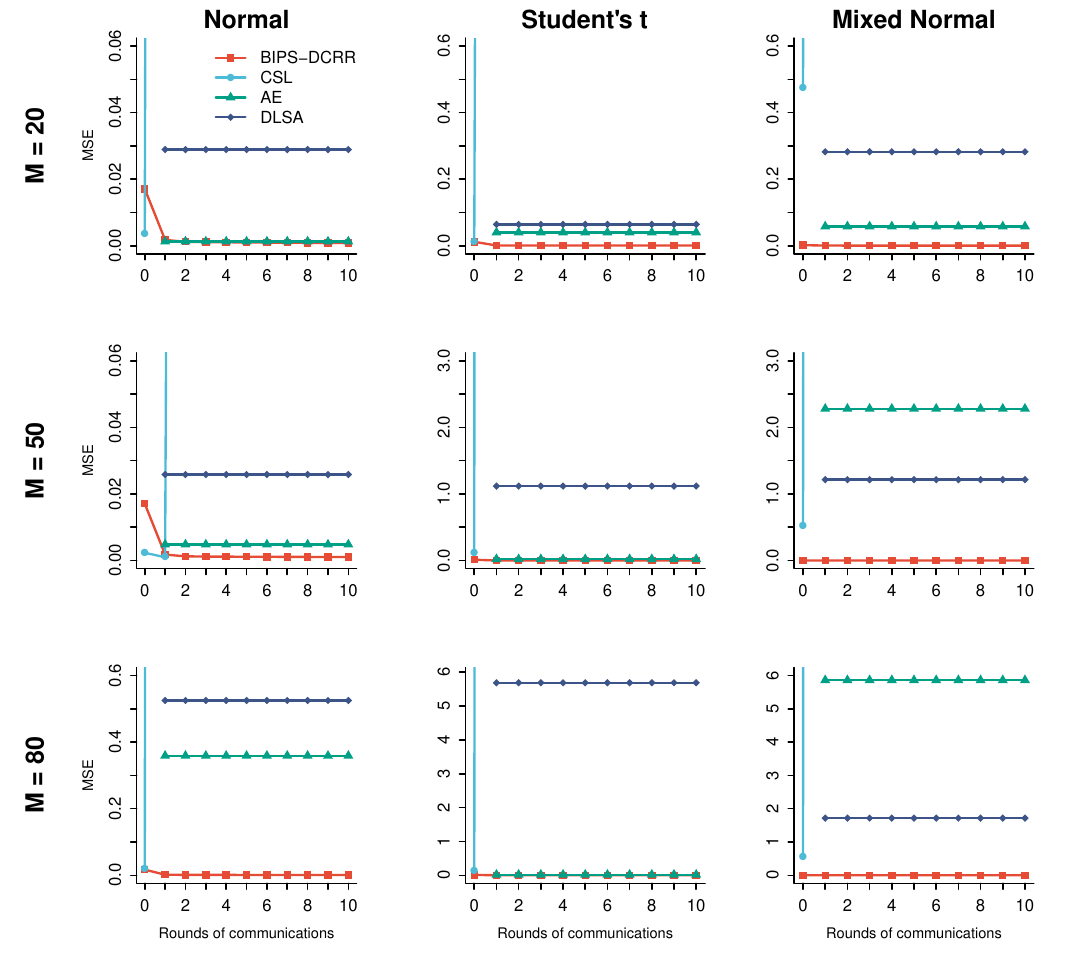} \\
\caption{Results of MSEs for four distributed estimation methods under Scheme 2, where we fix $n=50,000, r_0=200, r=1,000$, $\varrho=0.35$, $h=1$ and $\mathbf{X}\sim \mathrm{AR}(0.5)$.  (Normal): $\epsilon\sim \mathcal{N}(0,1)$, (Student's t): $\epsilon\sim \sqrt{2}t(4)$, and (Mixed Normal): $\epsilon\sim 0.95\mathcal{N}(0,1)+0.05\mathcal{N}(0,50)$. }
\label{fig:506}
\end{figure}

The data distributed on machines are in a random order in Scheme 1, while those in Scheme 2 follow a deterministic trend. Under both schemes, we  set $r=1000$, $\varrho=0.35$ and $h=1$ in BIPS-DCRR. Figures \ref{fig:505} and \ref{fig:506} provide the results for the four distributed methods above under Scheme 1 and Scheme 2, respectively. We can see that overall, under three types of error distributions, different numbers of machines, and two data distributed storage schemes, our proposed BIPS-DCRR method can achieve the convergence of the minimum MSE after at most 2 to 3 rounds of communication. Our proposed method performed better than other competing methods, especially under the mixed normal error distribution and under Scheme 2. The CSL method actually does not work under Scheme 2 since it can not achieve convergence. These results indicate that although existing methods work reasonably well when the error follows a normal distribution under random schemes, our method has a significant advantage in handling distributed data with some non-random storage schemes. More simulation results are available in Appendix B.4 of supplementary file.

\section{Real Data Analysis}
\begin{figure} %[htp]
\centering
    \includegraphics[height=8cm, width=14cm]{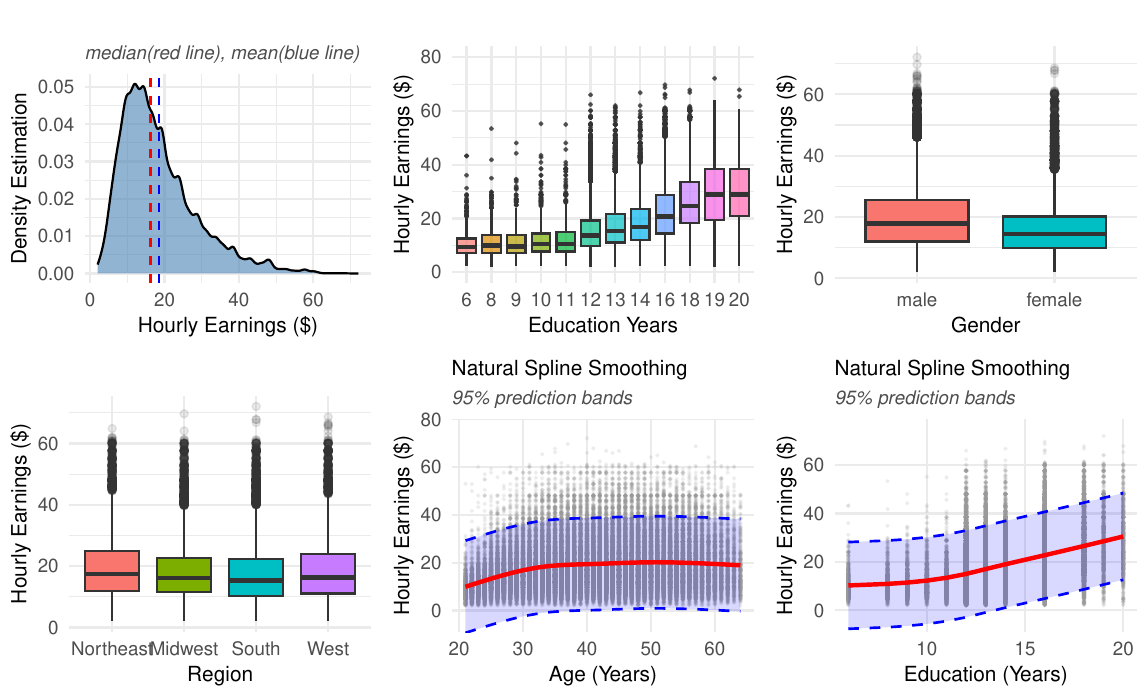} \\
\caption{Descriptive analysis plots for the CPSSW8 dataset. }
\label{fig:601a}
\end{figure}

%In this section, we apply our approach to the physicochemical properties of protein tertiary structure dataset (Dua and Graff, 2017), which contains 45,730 observations and the response variable is the size of the residue ranging from 0 to 21 Armstrong.
We apply our approaches to the CPSSW8 dataset, which is accessible via R package {\tt AER}. This dataset is extracted from the March U.S. Current Population Survey (CPS) administered by the Bureau of Labor Statistics and comprises a sample of $N = 61,395$ full-time individuals aged 21-64. The primary goal is to predict the average hourly earnings conditional upon a vector of covariates: education (years of completed schooling), age, gender, and geographic region. We carry out a descriptive analysis and the obtained results are displayed in Figure \ref{fig:601a}, revealing that the hourly earnings exhibits severe right-skewness and heavy-tailed distribution, wherein the sample mean is substantially inflated above the median by high-income outliers. Besides, % the distributions of covariates display that educational attainment exhibits pronounced mass points at discrete institutional milestones (e.g., 12 and 16 years), which correspond to expanding conditional wage variances, and systematic demographic and spatial heterogeneity is also manifest across gender and regional zones. 
the empirical distribution of education exhibits distinct mass points corresponding to high school (12 years) and bachelor’s degrees (16 years). The highly non-normal distribution displays the skewed empirical density of the response variable alongside conditional boxplots across education, gender, and regional categories. The bottom-right plots show nonlinear relations between the hourly earnings and education and age, separately. Afterwards, we first scale two numerical covariates (age and education) onto the interval $[0,1]$ with the transformation $(Z-\min_k Z_k)/(\max_{k}Z_k - \min_k Z_k)$, and apply the natural cubic spline with four degrees of freedom to each. We convert the categorical covariates (gender and region) into dummy variables, thereby generating a total of 12 predictors. The response variable is taken as the logarithmic hourly earnings. A standard regression analysis detects a total of 451 potential outlying data points. After removing the outliers, we have 60,944 remaining samples entering the downstream analysis.  

To evaluate the predictive performance of our subsampling methods with two existing methods, we use three metrics: mean absolute prediction errors ($\mathrm{MAPE}= \sum_{i \in \mathcal{T}} |y_i-\widehat{y}_i|/|\mathcal{T}|$), mean squared prediction errors ($\mathrm{MSPE}= \sum_{i \in \mathcal{T}} (y_i-\widehat{y}_i)^2/|\mathcal{T}|$), and the out-of-sample  $R^2$ ($R_{\text{oos}}^2=1-\sum_{i\in \mathcal{T} } (y_i-\widehat{y}_i)^2/\sum_{i\in \mathcal{T} } (y_i-\bar{y})^2$), where $ \mathcal{T}$ denotes the index set of test data, $\widehat{y}_i$ is the predicted response at the $i$th test data point and $\bar{y}$ is the mean of the response in train set. The whole dataset is randomly divided into a training set of $N_{\text{train}} = 50,000$ observations and a test set of the remaining $N_{\text{test}} = 10,944$ observations. We report the mean values of three metrics across 100 random partitions of the data. We fix the parameter  $\varrho = 0.35$ and the pilot sample size $r_0=200$. We compare distinct subsampling methods under three values of the Poisson subsample size $r$. The results are summarized in Table \ref{tab:601a}, showing that the proposed BIPS-CRR(A-opt) and BIPS-CRR(L-opt) are two best methods and all the BIPS methods outperform two benchmark methods (PS-QR and PS-LM).  

\begin{table}[htbp]
\centering
\footnotesize
\caption{Predictive performance of ten subsampling methods on the CPSSW8 dataset. The bold numbers indicate the top two results. }
\label{tab:601a}
\begin{tabular}{lccccccccc}
\toprule
 & \multicolumn{3}{c}{$r=1000$} & \multicolumn{3}{c}{$r=2000$} & \multicolumn{3}{c}{$r=3000$} \\   \cmidrule(r){2-4} \cmidrule(r){5-7}  \cmidrule(r){8-10}
Methods & MAPE & MSPE & $R^2_{\text{oos}}$ & MAPE & MSPE & $R^2_{\text{oos}}$ & MAPE & MSPE & $R^2_{\text{oos}}$ \\
\hline
BIPS-RR(L-opt)  & 0.35840 & 0.20628 & 0.29144 & 0.35747 & 0.20522 & 0.29506 & \textbf{0.35710} & 0.20483 & 0.29642 \\
BIPS-RR(A-opt)  & 0.35858 & 0.20647 & 0.29078 & 0.35747 & 0.20528 & 0.29489 & 0.35713 & 0.20490 & 0.29618 \\
BIPS-RR(lev)    & 0.35889 & 0.20673 & 0.28988 & 0.35765 & 0.20541 & 0.29443 & 0.35727 & 0.20500 & 0.29582 \\
BIPS-RR(unif)   & 0.35874 & 0.20658 & 0.29041 & 0.35770 & 0.20548 & 0.29417 & 0.35725 & 0.20502 & 0.29575 \\
BIPS-CRR(L-opt) & \textbf{0.35838} & \textbf{0.20617} & \textbf{0.29182} & \textbf{0.35738} & \textbf{0.20510} & \textbf{0.29549} & \textbf{0.35710} & \textbf{0.20480} & \textbf{0.29652} \\
BIPS-CRR(A-opt) & \textbf{0.35821} & \textbf{0.20603} & \textbf{0.29229} & \textbf{0.35729} & \textbf{0.20502} & \textbf{0.29576} & \textbf{0.35700} & \textbf{0.20472} & \textbf{0.29680} \\
BIPS-CRR(lev)   & 0.35875 & 0.20660 & 0.29034 & 0.35755 & 0.20531 & 0.29478 & 0.35719 & 0.20495 & 0.29600 \\
BIPS-CRR(unif)  & 0.35876 & 0.20664 & 0.29019 & 0.35760 & 0.20539 & 0.29451 & 0.35722 & 0.20494 & 0.29604 \\
PS-QR      & 0.35987 & 0.20800 & 0.28552 & 0.35812 & 0.20610 & 0.29227 & 0.35762 & 0.20551 & 0.29407 \\
PS-LM           & 0.35942 & 0.20748 & 0.28732 & 0.35800 & 0.20588 & 0.29282 & 0.35741 & 0.20526 & 0.29493 \\
\bottomrule
\end{tabular}
\end{table}

We finally investigate the predictive performance of proposed distributed method. We select $50,000$ observations as the training data, which are equally distributed over $M$ local machines, and the remaining observations as test data. We compare our BIPS-DCRR, the CSL method and the DLSA method, where the number of communications is fixed at 5 for BIPS-DCRR and CSL. We report the average of the MAPEs over 50 splits. The results are 0.3572 ($M = 25$) and 0.3569 ($M = 50$) for our BIPS-DCRR, which are smaller than 0.3575 ($M = 25$) and 0.3573 ($M = 50$) for CSL and 0.3587 ($M = 25$) and 0.3604 ($M = 50$) for DLSA.

\section{Discussion}
In this paper, we proposed a new subsampling framework, BIPS, for U-statistic empirical risk minimization and contributed two concrete estimation approaches (BIPS-CRR and BIPS-DCRR) under concoluted rank regression models with massive data. %The R code can be available via \url{https://github.com/statxxc/BIPS-DCRR}.
While our BIPS is primarily designed for CRR models, its underlying principles can be readily extended to other statistical models and machine learning tasks that involve pairwise loss functions. To illustrate this versatility, we provide several examples.
% Although our BIPS appears to be designed for CRR models, it can be readily extended to other statistical models and machine learning with pairwise loss functions. We provide a few examples for illustration. 
First, estimating the covariance $\sigma_{X,Y}$ of two random variables $X\in \mathbb{R}$ and $Y\in \mathbb{R}$ is one of fundamental tasks in risk modeling. The full estimator of $\sigma_{X,Y}$ can be casted as the minimizer of the pairwise loss function $\frac{1}{n(n-1)} \sum_{i=1}^{n}\sum_{j\neq i} [\frac{1}{2}(X_i-X_j)(Y_i-Y_j)-\sigma_{X,Y}]^2$ \citep{BC2018}. Thus, one can construct a BIPS estimator of $\sigma_{X,Y}$ as $ \argmin_{\sigma_{X,Y}} \frac{1}{n(n-1)} \sum_{i=1}^{n}\sum_{j\neq i} W_{i}W_{j} [\frac{1}{2}(X_i-X_j)(Y_i-Y_j)-\sigma_{X,Y}]^2$. 
Second, the pairwise learning including metric learning, ranking and online learning is an important area in machine learning. One of supervised distance metric learning problems is to solve the following pairwise learning problem \citep{JWZ2009}: $\min_{D\in \mathbb{R}_{+}^{p\times p}} \frac{1}{n(n-1)}\sum_{i=1}^{n}\sum_{j\neq i} \phi(Y_iY_j (1-d_{D}(\mathbf{X}_i, \mathbf{X}_j)))$, where $\phi$ is any loss function in maximum margin classifier, $Y_i\in \{-1, 1\}$ stands for the class label, $\mathbf{X}_i\in \mathbb{R}^{p}$ denotes the vector of predictors for the $i$-th data point, $d_{D}(\mathbf{X}_i, \mathbf{X}_j)=(\mathbf{X}_i-\mathbf{X}_j)^{\top}D(\mathbf{X}_i-\mathbf{X}_j)$ means the Mahalanobis distance between two points $\mathbf{X}_i$ and $\mathbf{X}_j$, and $D\in \mathbb{R}_{+}^{p\times p}$ is the distance metric to be learned. Naturally, a BIPS estimator of $D$ can be formulated as $\argmin_{D\in \mathbb{R}_{+}^{p\times p}} \frac{1}{n(n-1)}\sum_{i=1}^{n}\sum_{j\neq i} W_iW_j\phi(Y_iY_j (1-d_{D}(\mathbf{X}_i, \mathbf{X}_j)))$. In ranking learning, \cite{CLV2008} considered the following ranking problem. Let $\gamma$ be a ranking rule, namely, a function that maps $\mathbb{R}^{p}\times \mathbb{R}^{p} \to \{-1, 1\}$. Given data $\mathcal{D}_n=\{(\mathbf{X}_i, Y_i)\}_{i=1}^{n}$, where $\mathbf{X}_i\in \mathbb{R}^{p}$ is a random object and $Y_i\in  \mathbb{R}$ is the corresponding real-valued label. Note that $\mathbf{X}_i$ is said to be ``better" than $\mathbf{X}_j$ if $Y_i> Y_j$, and $\gamma(\mathbf{X}_i,\mathbf{X}_2)=1$ ($-1$) means that the rule ranks $\mathbf{X}_i$ higher (lower) than $\mathbf{X}_j$.  Learning the rule $\gamma$ over a class of ranking rules, $\Gamma$, amounts to the following minimization: $\min_{\gamma \in \Gamma} \frac{1}{n(n-1)}\sum_{i=1}^{n}\sum_{j\neq i} I(Z_{i,j}\cdot \gamma(\mathbf{X}_i, \mathbf{X}_j)<0)$, where $Z_{i,j}=(Y_i-Y_j)/2$. Thus, a BIPS estimator of the rule $\gamma$ can be computed as $\argmin_{\gamma \in \Gamma} \frac{1}{n(n-1)}\sum_{i=1}^{n}\sum_{j\neq i} W_iW_jI(Z_{i,j}\cdot \gamma(\mathbf{X}_i, \mathbf{X}_j)<0)$. However, for all these BISP problems, deriving the optimal subsampling probabilities in the weights $\{W_i\}_{i=1}^{n}$ remains a theoretical challenge and falls outside the scope of this paper.

Moreover, one could make serval extensions further. (i) Extending the proposed BIPS-DCRR to high-dimensional sparse massive data by introducing regularization terms in the objective function is an interesting task. It is straightforward to implement such penalized estimators for BIPS-CRR or BIPS-RR. However, more theoretical investigations are needed to support the empirical applications. (ii) Decentralized distributed estimation for CRR could further improve the computational performance when central machine is not available or can be easily conquered by attackers. (iii) Integrating data privacy-preserving mechanism such as differential privacy to subsampling and distributed estimation is still under developed.

%%%%%%%%%%%%%%%%%%%%%%%%%%%%%%%%%%%%%%%%%%%%%%%%%%%%%%%%%%%%%%%%%%%%%%%%%%%%%%%%%%%%%%%%%%%%%%%%%%%%%%%%%%%%%%%%%%%%%%%%%%%%


\begin{thebibliography}{99}
\bibitem[\protect\citeauthoryear{Ai et al.}{2021a}]{Aetal2021a} Ai, M., Wang, F., Yu, J., and Zhang, H. (2021a) Optimal subsampling for large-scale quantile regression. \textit{Journal of Complexity} \textbf{62}, 101512.

\bibitem[\protect\citeauthoryear{Ai et al.}{2021b}]{Aetal2021b} Ai, M., Yu, J., Zhang, H. and Wang, H. (2021b) Optimal subsampling algorithms for big data regressions. \textit{Statistica Sinica} \textbf{31}, 749-772.

% \bibitem[\protect\citeauthoryear{Bickel, Gotze and van Zwet}{1997}]{BGZ1997} Bickel, P., Gotze, F. and van Zwet, W. (1997) Resampling fewer than n observations: gains, losses, and remedies for losses. \textit{Statistica Sinica}, \textbf{7}, 1-31

\bibitem[\protect\citeauthoryear{Bose and Chatterjee}{2018}]{BC2018}  Bose, A. and Chatterjee, S. (2018) \textit{$U$-Statistics, $M_m$-Estimators and Resampling}. Springer Nature Singapore Pte Ltd.

%\bibitem[\protect\citeauthoryear{Chao et al.}{2024}]{Cetal2024} Chao, Y., Huang, L., Ma, X. and Sun, J. (2024) Optimal subsampling for modal regression in massive data. \textit{Metrika} \textbf{87}(4), 379-409.


\bibitem[\protect\citeauthoryear{Cl\'{e}mençon et al.}{2008}] {CLV2008} Cl\'{e}mençon, S., Lugosi, G. and Vayatis, N. (2008) Ranking and empirical minimization of U-statistics. \textit{The Annals of Statistics}, \textbf{36}(2), 844-874.

\bibitem[\protect\citeauthoryear{Chao et al.}{2024}]{CMZ2024} Chao, Y., Ma, X. and Zhu, B. (2024) Distributed optimal subsampling for quantile regression with massive data. \textit{Journal of Statistical Planning and Inference}, \textbf{233}, 106186.

\bibitem[\protect\citeauthoryear{Chao et al.}{2026}]{CXZ2026} Chao, Y., Xia, X. and Zhong, W. (2026) Communication-efficient pilot estimation for non-randomly distributed data in diverging dimensions. \textit{Journal of Computational and Graphical Statistics}, \textbf{35}(1), 155-172.


\bibitem[\protect\citeauthoryear{Chen et al.}{2020}]{Cetal2020}  Chen, X., Liu, W., Mao, X. and Yang, Z. (2020) Distributed high-dimensional regression under a quantile loss function. \textit{Journal of Machine Learning Research}, \textbf{21}, 1-43

\bibitem[\protect\citeauthoryear{Cook et al.}{2019}]{Cetal2019} Cook, C. E., Lopez, R., Stroe, O., Cochrane, G., et al.  (2019). The European bioinformatics institute in 2018: tools, infrastructure and training. \textit{Nucleic Acids Research}, \textbf{47}(1), 15-22.  %Brooksbank, C., Birney, E. and Apweiler, R.


\bibitem[\protect\citeauthoryear{Fan et al.}{2019}]{Fetal2019} Fan, J. Q., Wang, D., Wang, K. Z. and Zhu, Z. W. (2019) Distributed estimation of principal eigenspaces. \textit{The Annals of Statistics}, \textbf{47}(6), 3009-3031


\bibitem[\protect\citeauthoryear{Fan et al.}{2023}]{FGW2023} Fan, J., Guo, Y. and Wang, K. (2023) Communication-efficient accurate statistical estimation. \textit{Journal of the American Statistical Association}, \textbf{118}, 1000-1010


\bibitem[\protect\citeauthoryear{Gao et al.}{2022}]{Getal2022} Gao, Y., Liu, W., Wang, H., Wang, X., Yan, Y. and Zhang, R. (2022). A review of distributed statistical inference. \textit{Statistical Theory and Related Fields}, \textbf{6}(2), 89-99.



%\bibitem[\protect\citeauthoryear{Chung}{2001}]{C2001} Chung, K.L. (2001). \textit{A Course in Probability Theory}. Academic Press.




\bibitem[\protect\citeauthoryear{He et al.}{2024}]{Hetal2024} He, L., Li, W., Song, D. and Yang, M.S. (2024) A systematic view of information-based optimal subdata selection: algorithm development, performance evaluation, and application in financial data. \textit{Statistica Sinica}, \textbf{34}, 611-636.


\bibitem[\protect\citeauthoryear{Hesterberg}{1995}]{H1995}  Hesterberg, T. (1995) Weighted average importance sampling and defensive mixture distributions. \textit{Technometrics}, \textbf{37}(2), 185-194.

\bibitem[\protect\citeauthoryear{He and Xia}{2025}]{HX2025}  He, S. and Xia, X. (2025) Random perturbation subsampling for rank regression with massive data. \textit{Statistics and Computing}, \textbf{35}, 14.

% \bibitem[\protect\citeauthoryear{He et al.}{2021}]{Hetal2021}  He, X., Pan, X., Tan, K. M. and Zhou, W. X. (2021) Smoothed quantile regression with large-scale inference. \textit{Journal of Econometrics}, \textbf{232}, 367-388.

\bibitem[\protect\citeauthoryear{Hettmansperger and McKean}{2010}]{HM2010}  Hettmansperger, T. P. and McKean, J. W. (2010) \textit{Robust Nonparametric Statistical Methods (Second Edition)}. Boca Raton, FL: CRC Press.

\bibitem[\protect\citeauthoryear{Huang and Huo}{2019}]{HH2019} Huang, C. and Huo, X. (2019). A distributed one-step estimator. \textit{Mathematical Programming}, \textbf{174}, 41-76.


\bibitem[\protect\citeauthoryear{Jaeckel}{1972}]{J1972} Jaeckel, L. A. (1972) Estimating regression coefficients by minimizing the dispersion of the residuals. \textit{The Annals of Mathematical Statistics}, \textbf{43}(5), 1449-1458.

\bibitem[\protect\citeauthoryear{Jin et al.}{2009}] {JWZ2009} Jin, R., Wang, S. and Zhou, Y. (2009) Regularized distance metric learning: theory and algorithm. In \textit{Proceedings of the 23rd International Conference on Neural Information Processing Systems}, pp862-870.
 
\bibitem[\protect\citeauthoryear{Jordan et al.}{2019}]{JLY2019}  Jordan, M.I., Lee, J.D. and Yang, Y. (2019). Communication-efficient distributed statistical inference. \textit{Journal of the American Statistical Association}, \textbf{114}(526), 668-681.

\bibitem[\protect\citeauthoryear{Joseph and Vakayil}{2022}]{JV2022} Joseph, V. R. and Vakayil, A. (2022)  Split: an optimal method for data splitting. \textit{Technometrics}, \textbf{64}, 166-176.


%\bibitem[\protect\citeauthoryear{Knight}{1998}]{K1998} Knight, K. (1998) Limiting distributions for L1 regression estimators under general conditions. \textit{The Annals of Statistics}, \textbf{26}, 755-770



% \bibitem[\protect\citeauthoryear{Koenker}{2005}]{K2005} Koenker, R. (2005). \textit{Quantile Regression}. Cambridge University Press, UK.

% \bibitem[\protect\citeauthoryear{Kar et al.}{2013}]{Ketal2013} Kar, P., Sriperumbudur, B., Jain, P. and Karnick, H. (2013) On the generalization ability of online learning algorithms for pairwise loss functions. In \textit{International Conference on Machine Learning}, 441-449.


\bibitem[\protect\citeauthoryear{Leng}{2010}]{L2010}  Leng, C. L. (2010) Variable selection and coefficient estimation via regularized rank regression. \textit{Statistica Sinica}, \textbf{20}, 167-181.


\bibitem[\protect\citeauthoryear{Li et al.}{2024}]{LXZ2024a}  Li, X., Xia, X. and Zhang, Z. (2024) Poisson subsampling-based estimation for growing-dimensional expectile regression in massive data. \textit{Statistics and Computing}, \textbf{34}, 133.



\bibitem[\protect\citeauthoryear{Lin and Li}{2023}]{LL2023} Lin, L. and Li, F. (2023) Global debiased DC estimations for biased estimators via pro forma regression. \textit{Test}, \textbf{32}(2), 726-758.


\bibitem[\protect\citeauthoryear{Liu et al.}{2025}]{LMT2025} Liu, W., Mao, X. and Tu, J. (2025) Communication-efficient distributed sparse learning with oracle property and geometric convergence. \textit{Journal of the American Statistical Association}, DOI:10.1080/01621459.2025.2479237 



% \bibitem[\protect\citeauthoryear{Li, Xia and Zhang}{2024b}]{LXZ2024b}  Li, X., Xia, X. and Zhang, Z. (2024b) Distributed subsampling for multiplicative regression. \textit{Statistics and Computing}, \textbf{34}(5), 1-20.


% \bibitem[\protect\citeauthoryear{Ma et al.}{2015}]{MMY2015} Ma, P., Mahoney, M. and Yu, B. (2015)  A statistical perspective on algorithmic leveraging. \textit{ Journal of Machine Learning Research} \textbf{16}(27), 861-919.







\bibitem[\protect\citeauthoryear{Ma et al.}{2015}]{MMY2015} Ma, P., Mahoney, M. and Yu, B. (2015) A statistical perspective on algorithmic leveraging. \textit{Journal of Machine Learning Research}, \textbf{16}(27), 861-911.

% \bibitem[\protect\citeauthoryear{McCullagh and Nelder}{1989}]{MN1989} McCullagh, P. and Nelder, J. A. (1989). \textit{Generalized Linear Models, 2nd ed}. Chapman and Hall, London.


\bibitem[\protect\citeauthoryear{Ren et al.}{2023}]{RZW2023}  Ren, M., Zhao, S. and Wang, M. (2023)  Optimal subsampling for least absolute relative error estimators with massive data.  \textit{ Journal of Complexity} \textbf{74}, 101694.


%\bibitem[\protect\citeauthoryear{Ren et al.}{2024}]{Retal2024} Ren, M., Zhao, S., Wang, M., and Zhu, X. (2024). Robust optimal subsampling based on weighted asymmetric least squares.\textit{ Statistical Papers}, \textbf{65}(4), 2221-2251.

%\bibitem[\protect\citeauthoryear{Ting and Brochu}{2018}]{TB2018} Ting, D. and Brochu, E. (2018) Optimal subsampling with influence functions, \textit{in Advances in Neural Information Processing Systems}, 31, Curran Associates, Inc., pp. 3654-3663.

%\bibitem[\protect\citeauthoryear{Serfling}{1980}]{S1980} Serfling, R. J. (1980). \textit{Approximation Theorems of Mathematical Statistics}. John Wiley \& Sons. Inc.


% \bibitem[\protect\citeauthoryear{Shao and Wang}{2022}]{SW2022} Shao, Y. and Wang, L. (2022). Optimal subsampling for composite quantile regression model in massive data. \textit{ Statistical Papers}, \textbf{63}(4), 1139-1161


\bibitem[\protect\citeauthoryear{Shao et al.}{2025}]{SWL2024} Shao, Y., Wang, L., and Lian, H. (2025). Optimal decorrelated score subsampling for high-dimensional generalized linear models under measurement constraints. \textit{Journal of Computational and Graphical Statistics}, \textbf{34}(2),  530-539

\bibitem[\protect\citeauthoryear{Su\'{a}rez et al.}{2021}]{SGH2021} Su\'{a}rez, J. L.,  García, S. and Herrera, F. (2021) A tutorial on distance metric learning: Mathematical foundations, algorithms, experimental analysis, prospects and challenges. \textit{Neurocomputing}, \textbf{425}, 300-322

%\bibitem[\protect\citeauthoryear{Shao, Wang and Lian}{2025}]{SWL2025} Shao, Y., Wang, L., and Lian, H. (2025). Optimal distributed subsampling under heterogeneity. \textit{Statistics and Computing}, \textbf{35}(2), 26.


\bibitem[\protect\citeauthoryear{Tan et al.}{2022}]{TBZ2022} Tan, K. M., Battey, H. and Zhou, W.X. (2022) Communication-constrained distributed quantile regression with optimal statistical guarantees. \textit{Journal of Machine Learning Research}, \textbf{23}, 12456-12516

\bibitem[\protect\citeauthoryear{Tan et al.}{2022}]{TWZ2022}  Tan, K. M., Wang, L., and Zhou, W.X. (2022). High-dimensional quantile regression: convolution smoothing and concave regularization. \textit{Journal of the Royal Statistical Society, Series B}, \textbf{84}, 205-233.





\bibitem[\protect\citeauthoryear{Wang}{2018}]{W2018} Wang, H. (2018) More efficient estimation for logistic regression with optimal subsamples.  \textit{Journal of Machine Learning and Research}, \textbf{20}(132), 1-59.


\bibitem[\protect\citeauthoryear{Wang and Li}{2009}]{WL2009}  Wang, L. and Li, R. (2009) Weighted Wilcoxon-type smoothly clipped absolute deviation method. \textit{Biometrics}, \textbf{65}, 564-571

\bibitem[\protect\citeauthoryear{Wang et al.}{2020}]{Wetal2020} Wang, L., Peng, B., Bradic, J., Li, R. and Wu, Y. (2020) A tuning-free robust and efficient approach to high-dimensional regression (with discussion). \textit{Journal of the American Statistical Association}, \textbf{115}, 1700-1714.



% \bibitem[\protect\citeauthoryear{Wang, Wu and Li}{2012}]{WWL2012}  Wang, L., Wu, Y. and Li, R. (2012) Quantile regression for analyzing heterogeneity in ultra-high dimension. \textit{Journal of the American Statistical Association}, \textbf{107}, 214-222.


\bibitem[\protect\citeauthoryear{Wang and Ma}{2021}]{WM2021} Wang, H. and Ma, Y. (2021). Optimal subsampling for quantile regression in big data. \textit{Biometrika}, \textbf{108}, 99-112.

\bibitem[\protect\citeauthoryear{Wang et al.}{2019}]{WYS2019} Wang, H., Yang, M. and Stufken, J. (2019) Information-based optimal subdata selection for big data linear regression. \textit{Journal of the American Statistical Association}, \textbf{114}(525), 393-405.

\bibitem[\protect\citeauthoryear{Wang et al.}{2019}]{Wetal2019} Wang, X., Yang, Z., Chen, X. and Liu, W. (2019) Distributed inference for linear support vector machine. \textit{Journal of Machine Learning Research},  \textbf{20}, 1–41

\bibitem[\protect\citeauthoryear{Wang et al.}{2018}]{WZM2018} Wang, H., Zhu, R. and Ma, P. (2018). Optimal subsampling for large sample logistic regression. \textit{Journal of the American Statistical Association}, \textbf{113}(522), 829-844.


\bibitem[\protect\citeauthoryear{Wang et al.}{2022}]{WZW2022} Wang, J., Zou, J. and Wang, H. (2022). Sampling with replacement vs Poisson sampling: A comparative study in optimal subsampling. \textit{IEEE Transactions on Information Theory}, \textbf{68}:10, 6605-6630.

\bibitem[\protect\citeauthoryear{Wen et al.}{2025}]{Wetal2025} Wen, W.,  Li, H., Wu, R., Wu, L. and Chen, H. (2025) Generalization analysis of adversarial pairwise learning, \textit{Neural Networks}, \textbf{183}, 106955

%\bibitem[\protect\citeauthoryear{Yao and Wang}{2021}]{YW2021} Yao, Y. and Wang, H. (2021). A review on optimal subsampling methods for massive datasets. \textit{Journal of Data Science}, \textbf{19}(1), 151-172.


\bibitem[\protect\citeauthoryear{Xu et al.}{2024}]{Xetal2024} Xu, W. L., Liu, J. M. and Lian, H. (2024) Distributed estimation of support vector machines for matrix data. \textit{IEEE Transactions on Neural Networks and Learning Systems}, \textbf{35}(5), 6643-6653



\bibitem[\protect\citeauthoryear{Yu et al.}{2024}]{YAY2024} Yu, J., Ai, M. and Ye, Z. (2024). A review on design inspired subsampling for big data. \textit{Statistical Papers}, \textbf{65}, 467-510.

\bibitem[\protect\citeauthoryear{Yu et al.}{2022}]{Yetal2022} Yu, J., Wang, H.Y., Ai, M. and Zhang, H. (2022) Optimal distributed subsampling for maximum quasi-likelihood estimators with massive data. \textit{Journal of the American Statistical Association}, \textbf{117}(537), 265-276.

\bibitem[\protect\citeauthoryear{Zhang et al.}{2013}]{ZDW2013} Zhang, Y. C., Duchi, J. C. and Wainwright, M. J. (2013) Communication-efficient algorithms for statistical optimization. \textit{Journal of Machine Learning Research}, \textbf{14}(11), 3321-3363.


% \bibitem[\protect\citeauthoryear{Zhang and Wang}{2021}]{ZW2021} Zhang, H. and Wang, H. (2021). Distributed subdata selection for big data via sampling-based approach. \textit{Computational Statistics \& Data Analysis}, \textbf{153}, 107072.

\bibitem[\protect\citeauthoryear{Zhou et al.}{2024}]{ZWZ2024}  Zhou, L., Wang, B. and Zou, H. (2024) Sparse convoluted rank regression in high dimensions. \textit{Journal of the American Statistical Association}, \textbf{119}(546), 1500-1512


\bibitem[\protect\citeauthoryear{Zhou et al.}{2023}]{Zetal2023} Zhou, Z., Yang, Z., Zhang, A. and Zhou, Y. (2023) Efficient model-free subsampling method for massive data. \textit{Technometrics}, \textbf{66}, 240-252.




\bibitem[\protect\citeauthoryear{Zhu et al.}{2021}]{ZLW2021}  Zhu, X., L, F. and Wang, H. (2021) Least-square approximation for a distributed system. \textit{Journal of Computational and Graphical Statistics}, \textbf{30}(4), 1004-1018.
\end{thebibliography}
\end{document}